\documentclass[11pt]{amsart}
\usepackage{mathrsfs}
\usepackage{color}
\usepackage{amssymb}
\usepackage{algorithmic}
\usepackage{algorithm}
\usepackage{rotating}
\usepackage{boldline}
\usepackage{amsmath}
\usepackage{bm}
\usepackage{graphicx}
\usepackage{verbatim}
\usepackage{slashbox}

\usepackage{float}

\usepackage{caption}
\usepackage{subcaption}
\usepackage{enumerate}
\begin{document}


\def \qn{d}
\def \qs{s}
\def \qkappa{\kappa}
\def \qmu{\mu}
\def \laplace{\Delta}
\def \qetauh{e^\tau}
\def \qbetauh{\bm{e}^\tau}
\def \qdbeh{\bm{e}_\delta}
\def \qdeh{e_\delta}
\def \qbeHh{\bm{e}_{H,\delta}}
\def \qeHh{e_{H,\delta}}
\def \qbeH{\bm{e}_H}
\def \qeH{e_H}
\def \qecalH{{F}}
\def \qbuH{\bm{u}_H}
\def \qbuh{\bm{u}_{\text{ref}}}
\def \qbu{\bm{u}}
\def \qbw{\bm{w}}
\def \qu{u}
\def \qbvarphi{\bm{\varphi}}
\def \qbv{\bm{v}}
\def \qv{v}
\def \qbvH{\bm{v}_H}
\def \qdbvh{\bm{v}_\delta}
\def \qdbuh{\qbu_\delta}
\def \qdph{p_\delta}
\def \qbvhtau{\bm{v}^{\tau}}
\def \qqhtau{q^{\tau}}
\def \qdtbuh{{\bm{u}}_\delta}
\def \qbuHh{\bm{u}_{H,\delta}}
\def \qbvHh{\bm{v}_{H,\delta}}
\def \qpH{p_H}
\def \qp{p}
\def \qq{q}
\def \qi{i}
\def \qj{j}
\def \qk{k}
\def \qr{r}
\def \qJscr{\mathscr{J}}
\def \qqH{q_H}
\def \qdqh{q_\delta}
\def \qph{p_{\text{ref}}}
\def \qdtph{{p}_\delta(\qbuH)}
\def \qepsilon{\epsilon}
\def \qpHh{p_{H,\delta}}
\def \qqHh{q_{H,\delta}}
\def \qdVscrh{\delta\mathscr{V}}
\def \qdWscrh{\delta\mathscr{W}}
\def \qtauVscrh{\mathscr{V}^{\tau}}
\def \qtauWscrh{\mathscr{W}^{\tau}}
\def \qbVscrH{\mathscr{V}_H}
\def \qWscrH{\mathscr{Q}_H}
\def \qWscrh{\mathscr{Q}_h}
\def \qWscrj{\mathscr{Q}_j}
\def \qTcalH{\mathcal{T}_h}
\def \qTcalj{\mathcal{T}_j}
\def \qTcalJ{\mathcal{T}_J}
\def \qVscrHh{\mathscr{V}_{H,\delta}}
\def \qVscrH{\mathscr{V}_{h}}
\def \qVscrj{\mathscr{V}_{j}}
\def \qVscrHd{\mathscr{V}_h^\partial}
\def \qWscrHh{\mathscr{Q}_{H,\delta}}
\def \qTH{{T}_h}
\def \qEcalH{\mathcal{F}_h}
\def \qEcalj{\mathcal{F}_j}
\def \qEcal{\mathcal{F}}
\def \qEH{N}
\def \qiEcal{{\mathcal{F}}_h}
\def \qiEcalH{{\mathcal{F}}^0_h}
\def \qiEcalj{{\mathcal{F}}^0_j}
\def \qniE{n_{\qiEcalH}}
\def \qdEcalH{{\mathcal{F}_h^\partial}}
\def \qdEcalj{{\mathcal{F}_j^\partial}}
\def \qdTH{\partial\qTH}
\def \qH{h}
\def \q0{\bm{0}}
\def \qbf{\bm{f}}
\def \qbfm{\bm{f}_m}
\def \qfs{f_s}
\def \qbg{\bm{g}}
\def \qg{g}
\def \qbubg{\qbu_{\qbg}}
\def \qbug{\qbu_{\qg}}
\def \qOmega{\Omega}
\def \qdOmega{\partial\Omega}
\def \qdbx{~d\bm{x}}
\def \qds{~ds}
\def \grad{\nabla}
\def \div{\nabla\cdot}
\def \qbn{\bm{n}}
\def \qbtau{\bm{\tau}}
\def \qalpha{\alpha}
\def \qFhmB{F_{m}}
\def \qFsB{F_s}
\def \qdbph{{p}_\delta(\qFhmB)}
\def \qdbbuh{{\bm{u}}_\delta(\qFhmB)}
\def \qC{C}
\def \qCe{C_{\iota}}
\def \qbx{\bm{x}}
\def \qx{x}
\def \qP{P}
\def \eps{{\varepsilon}}
\def \qbR{\mathbb{R}}
\def \qOcal{\mathcal{O}}
\def\dblLbrace{\{\hspace{-.45ex}\{}
\def\dblRbrace{\}\hspace{-.45ex}\}}

\definecolor{mid_gray}{rgb}{0.5,0.5,0.5}
\definecolor{light_gray}{rgb}{0.7,0.7,0.7}

\def\naught#1{#1^0}
\def\A{\mathcal A}
\def\B{\mathcal B}
\def\L{\mathcal L}
\def\I{\mathcal I}
\def\P{\mathcal P}
\def\X{\mathcal X}

\def\F{\mathbb F}
\def\T{\mathbb T}
\def\R{\mathbb R}
\def\Q{\mathbb Q}
\def\d{\partial}
\def\n{\mathbf n}
\def\grad{\nabla}
\def\div{\nabla\!\cdot\!}
\def\curl{\nabla\!\times\!}
\def\eps#1{\varepsilon(#1)}
\def\ntrace#1{#1\!\cdot\!\n}
\newcommand{\mvl}[1]{\{\!\!\{#1\}\!\!\}}             
\def\rf#1{\widehat{#1}}
 \def\Vh{V_}
\def\Hdiv{H^{\text{div}}}                            

\newcommand{\seminorm}[1]{\bigl|#1\bigr|}
\newcommand{\norm}[1]{\bigl\|#1\bigr\|}
\def\ipnorm#1{\norm{#1}_{1,h}}
\def\form(#1,#2){\bigl(#1,#2\bigr)}
\def\forme(#1,#2){\bigl<#1,#2\bigr>}

\newcommand{\jump}[1]{\mbox{$\left[\hspace{-2pt}\left[{#1}\right]\hspace{-2pt}\right]$}}

\newcommand{\enorm}[1]{\,|\!|\!| {{#1}} |\!|\!|}

\renewcommand{\(}{\left(}
\renewcommand{\)}{\right)}
\renewcommand{\[}{\left[}
\renewcommand{\]}{\right]}
\newcommand{\<}{\left<}
\renewcommand{\>}{\right>}
\newcommand{\tnorm}[2]{\interleave #1 \interleave_{#2}}
\newcommand{\dnorm}[2]{\left\| #1 \right\|_{#2}}
\newcommand{\snorm}[2]{\left| #1 \right|_{#2}}

\renewcommand{\dfrac}[2]{\displaystyle{\frac{#1}{#2}}}
\newcommand{\dint}[2]{\displaystyle{\int_{#1}^{#2}}}
\newcommand{\dcup}[2]{\displaystyle{\bigcup_{#1}^{#2}}}
\newcommand{\dsum}[2]{\displaystyle{\sum_{#1}^{#2}}}
\newcommand{\dmax}[1]{\displaystyle{\max_{#1}}}
\newcommand{\dmin}[1]{\displaystyle{\min_{#1}}}
\newcommand{\dsup}[1]{\displaystyle{\sup_{#1}}}
\newcommand{\dist}{\text{dist}}
\newcommand{\ngrad}[1]{\dfrac{\partial #1}{\partial\qbn}}
\newcommand{\average}[1]{\dblLbrace \eps(#1) \dblRbrace}

\newcommand{\qaH}[2]{a_H\( #1, #2\)}
\newcommand{\qb}[2]{b\( #1, #2\)}
\newcommand{\qah}[2]{a_{h}\( #1, #2\)}
\newcommand{\qaHh}[2]{a_{H,h}\( #1, #2\)}
\newcommand{\qaD}[2]{a^D\( #1, #2\)}
\newcommand{\qahB}[2]{a\( #1, #2\)}
\newcommand{\qaohB}[2]{a_{\qTH}^B\( #1, #2\)}
\newcommand{\qathB}[2]{\widetilde{a}\( #1, #2\)}
\newcommand{\qatoB}[2]{\mathring{\widetilde{a}}^B\( #1, #2\)}
\newcommand{\qatD}[2]{\widetilde{a}^D\( #1, #2\)}

\definecolor{gray}{rgb}{0.6,0.6,0.6}
\newcommand{\gray}[1]{\textcolor{gray}{#1}}
\newcommand{\red}[1]{\textcolor{red}{#1}}
\newcommand{\green}[1]{\textcolor{green}{#1}}
\newcommand{\blue}[1]{\textcolor{blue}{#1}}
\newcommand{\magenta}[1]{\textcolor{magenta}{#1}}

\newcounter{dummy} \numberwithin{dummy}{section}

\newtheorem{theorem}[dummy]{Theorem}
\newtheorem{lemma}[dummy]{Lemma}
\newtheorem{proposition}[dummy]{Proposition}
\newtheorem{corollary}[dummy]{Corollary}
\newtheorem{remark}[dummy]{Remark}
\newtheorem{example}[dummy]{Example}

\newlength{\figurewidth}
\setlength{\figurewidth}{\textwidth}
\newlength{\figurewidththree}
\setlength{\figurewidththree}{\textwidth}



\def\AAA{\bm A}
\def\B{\bm B}
\def\R{\bm R}

\title[MG for Darcy/Brinkman equations]{Geometric Multigrid for Darcy and Brinkman models of flows in 
highly heterogeneous porous media:\\ A numerical study}

\author{Guido Kanschat} 

\address{Interdisziplin\"ares Zentrum f\"ur Wissenschaftliches Rechnen (IWR),
Universit\"at Heidelberg, Im Neuenheimer Feld 368, 69120 Heidelberg, 
Germany ({\tt kanschat@uni-heidelberg.de}).}

\author{Raytcho Lazarov}

\address{Department of Mathematics, Texas A\&M University, College Station, Texas, 77843-3368, USA  \\
({\tt lazarov@math.tamu.edu}).}

\author{Youli Mao}

\address{Department of Geology \& Geophysics, Texas A\&M University, College Station, Texas, 77843-3368, USA \\
({\tt youlimao@tamu.edu}).}

\date{started April 20, 2013, today is \today}


\begin{abstract}
  We apply geometric multigrid methods for the finite element approximation of flow problems governed by Darcy and Brinkman  systems used in modeling highly heterogeneous porous media. The   method is based on divergence-conforming discontinuous Galerkin methods  and overlapping, patch based domain decomposition smoothers. We show in benchmark experiments that the method is robust with respect to  mesh size and contrast of permeability for highly heterogeneous media.
\end{abstract}


\maketitle

\thispagestyle{plain}

\section{Introduction}\label{sec:intro}
This paper is devoted to development and testing of discretization and
multilevel solution algorithms for a unified approach in simulations
involving Darcy and Brinkman
models of flows in highly
heterogeneous porous media.
Many processes in engineering, geophysics, and hydrology involve such
flows. They are modeled by systems of partial differential equations
that are similar to those used in heat transfer, diffusion,
filtration, and other industrial processes.  The common characteristic
of these diverse areas is that media may exhibit heterogeneities over
a wide range of length-scales.  Depending on the goals and the
particular applications the solution of the corresponding mathematical
problem might be sought at various scales.

For example, if in reservoir modeling we are interested in the global
pressure drop for a given flow rate (and no fine scale details of the
flow are important) then the problem is formulated and solved on a
{\it field (large)} scale with some average reservoir characteristics.
In this case we use some upscaling procedure or homogenized
mathematical model.  However, often one may need to have information
about some fine scale details. In such cases one can use two-scale or
multiscale methods which are capable of enriching the global
coarse-scale solution with fine-scale features. And finally, when all
fine-scale features are needed one should use a detailed model which
uses all available fine-scale information.  Since the microstructure
influences physical macro properties of cellular materials (e.g.
permeability, acoustic or thermal properties, stiffness, etc), the
solutions of such fine-scale models are often used to calculate these
properties.

Such applications have motivated our study of numerical methods and
algorithms for simulation of fluid flows in highly heterogeneous
porous media with low volume fraction of the solid matrix resulting in high
porosity. At pore level, the Reynolds number is small due to the
small reference length. Therefore, flows of incompressible fluids can
be modeled by Stokes' equations.  On a field-scale, fluid flows in
porous media have been modeled mainly by mass conservation equation
and by Darcy's law $ \nabla \qp = - \qmu K^{-1} \qbu$ between the
macroscopic pressure $\qp$ and velocity $\qbu$, which we
write in the form  $ \nabla \qp = - \kappa(x) \qbu$. Here $K(x)$ the media
permeability, $\qmu$ is the fluid viscosity, and $\kappa(x)$ is the scaled inverse permeability.

\begin{figure}[tp]
  \begin{subfigure}[H]{0.45\textwidth}
    \includegraphics[width=5.2cm, height=4.5cm] 
    {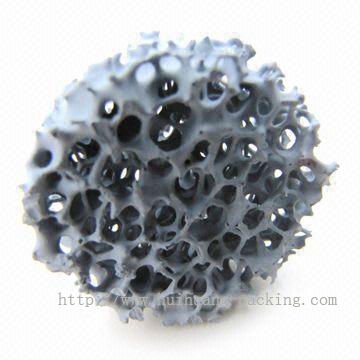} \hspace{1cm}
    \caption{CT scan of open foam}   
  \end{subfigure}
  \begin{subfigure}[H]{0.45\textwidth}
    \includegraphics[height = 4.5cm, width=7cm] 
    {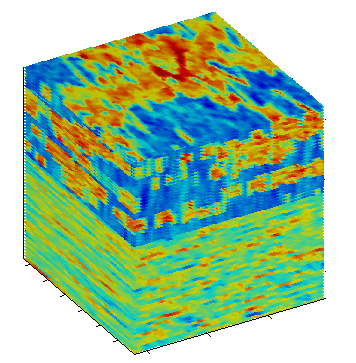}
    \caption{Permeability of SPE10 benchmark}   
  \end{subfigure}
  \caption{CT scans of highly porous materials on a micro-scale}
  \label{fig:introduction}
\end{figure}

Many porous media are characterized by very low solid volume fraction
and thus high porosity, e.g. fractured or vuggy reservoirs, mineral wool,
and industrial foams, (cf. Figure~\ref{fig:introduction} (A)). For
such media the porosity could be as high as 95 -- 98 \%.  For such
highly porous media Darcy's law often does not give good agreement
with the experimental data. In order to reduce the deviations between
the measurements for flows in highly porous media and the Darcy-based
predictions, Brinkman in \cite{Bri47_1} introduced a new
phenomenological relation between the velocity and the pressure
gradient $\nabla p= - 
\kappa(x) \qbu + \mu \Delta \qbu$, (see,
\cite[page 94]{Kaviany91}).  Together with conservation of mass, which
in the absence of any mass sources or sinks is expressed by
$ \nabla \cdot \qbu =0$, Brinkman's relations and proper boundary
conditions form a closed mathematical model \eqref{eq:Brinkman}.  An
important characteristic is the contrast of the media
$\boldsymbol{\kappa} $, defined as the ratio between the highest and
lowest values of the permeability,
$\boldsymbol{\kappa} = \max_{x \in \Omega} K(x)/  \min_{x \in \Omega}  K(x)$. 
The problems we consider in this paper involve $K(x)$ varying by
orders of magnitude on small length scales.

Darcy's and Brinkman's equations were introduced as phenomenological
macroscopic equations without direct link to underlying microscopic
structure of the media.  Advances in homogenization theory made it
possible to rigorously derive them from Stokes' equations for periodic
media, see e.g.\ \cite{All91}. As concluded in
\cite[pp.~266--273]{All91}, there are three different limits depending
on the size of the periodically arranged obstacles, which respectively
lead to Darcy's, Brinkman's, and Stokes' equations as macroscopic,
i.e., homogenized, relation.  There is still ongoing discussion about
the validity and the applications of Brinkman equations as a model of
flows in porous media at higher values of the porosity, see,
e.g. \cite{Auriault_2009,Durlofsky_1987}.  Nevertheless, Brinkman's
system of equations is  a convenient model of 
flows in highly heterogeneous porous media with random
distribution of inclusions, obstacles, channels, layers and other
geophysical features.

The Darcy/Brinkman model has been used also in the framework of 
fictitious domain methods that allows to treat in a unified
way flows in porous media, \cite{Konovalov79}, time dependent
incompressible viscous flows,
\cite{Bugrov_Smagulov,Glowinski94,KAPC00}, and transient compressible
viscous flows, \cite{Vabishchevich_FRM}.
A rigorous analysis of Brinkman's system from the point of view of the
fictitious domain method was carried out in \cite{Ang99,KAPC00}. As a
result, Stokes equations in a complicated domain (flow around many
obstacles, an obstacle with complicated topology,  or domains with void or caverns, e.g.\
Figure~\ref{fig:introduction} (A)), are replaced by Brinkman's
equations in a simpler domain but with highly varying 
$\kappa(x)$. In such models the values of $\kappa(x)$ 
in the obstacle/void plays a role of a penalty
parameter that is directly related to the contrast
$\boldsymbol{\kappa}$ so that large variations in $\kappa(x)$ lead to
ill-conditioning of the corresponding discrete problem.

In this paper we consider a unified approach to approximation of the
Brinkman/Darcy flow equations by $H^{\text{div}}$-conforming
Raviart-Thomas mixed finite elements
(cf.~\cite{CockburnKanschatSchoetzau07} for Navier-Stokes equations
and~\cite{GiraultKanschatRiviere14} for Darcy-Stokes coupling).  Close
to our research, discontinuous Galerkin (DG) FEM with 
$H^{\text{div}}$-conforming elements has
also been applied to numerically solve the Brinkman system in
~\cite{Konno_2011}, where a priori and a posteriori error analysis in
both Darcy and Brinkman limits are performed.  This discretization is
applied on fine meshes resolving mesoscale heterogeneities of the
media.  That is, we assume a homogenized microscale structure of
material below grid resolution and changes of this microstructure are
resolved.  For example, certain media are discretized with $256^d$ or
$512^d$ ($d=2,3$ is the space dimension) voxels. Discretization on
such grids results in very large algebraic saddle point problems which
are ill-conditioned due to both, the small mesh size and high contrast
$\boldsymbol{\kappa}$.

When the jumps in the permeability are aligned with a coarse mesh,
techniques based on domain decomposition, geometric multigrid or
multilevel methods could lead to efficient and optimal methods for
Darcy's model,
e.~g.~\cite{Kraus2009robust,powell2005parameter,powell2003optimal,spillane2014}.
A multigrid preconditioner for solving algebraic systems resulting from finite element 
approximation in $H^{\text{div}}$ was proposed 
by Arnold, Falk, and Winther in \cite{Arnold1997preconditioning,Arnold2000MG}. 
The analysis of the preconditioner relies on the regularity of the solution.
In various practical computations, it was demonstrated that this preconditioner is robust
with respect to the contrast $\boldsymbol{\kappa}$ for Darcy's model,
e.g. \cite[Tables 4,5]{powell2005parameter} and \cite[Tables 2.3 --
2.7]{powell2003optimal}. However, when the jumps in the permeability
are not aligned with the coarse mesh the multigrid (or multilevel)
method performance deteriorates when the contrast gets larger,
e.g. \cite[Table 7.10]{Kraus2009robust}.  Robust with respect to the contrast
$\boldsymbol{\kappa}$
multilevel preconditioners  based on additive Schur complement approximation
were recently proposed and experimentally studied in
\cite{Kraus_2015} for the Darcy system. 

The class of problems we consider in this paper is characterized by 
large variations of the permeability 
on a very fine scale (high frequency) which results in low regularity of the
solution. Also, the permeability is given on a fine scale and
there is no practical way to split the domain into a fixed number of
subdomains where the permeability is smooth or constant.  To solve
such system we use a multigrid preconditioner consistent with the
divergence free subspace and employing 
a smoother that uses patches around grid vertexes and based on the idea of 
Arnold, Falk, and Winther~\cite{Arnold1997preconditioning,Arnold2000MG}.
As shown recently in
\cite{KanschatMao_JNM,KanschatMao14}, this preconditioner is optimal
for divergence-conforming DG approximations
of Stokes' equations.  The goal of this paper is to study numerically
the performance of this multigrid approach to solving the
Darcy/Brinkman system for high contrast and high frequency porous
media.

 Preconditioning of DG FEM for Brinkman equations was
considered in \cite{efendiev2012robust}.  It was proven theoretically
and confirmed experimentally that in two space dimensions the proposed
preconditioner based on domain decomposition technique (or two-level
method) that involves solution of some local spectral problems is
optimal with respect to both, the mesh-size and the
contrast. Multilevel generalization of the same idea applied to
anisotropic problems was done in \cite{Willems_2013} and extended to
abstract symmetric and positive definite forms in \cite{Willems_2014}.
The numerical experiments in \cite{Willems_2013,Willems_2014} show
that the proposed method is robust with respect to both the mesh-size
and media contrast. However, the theory depends on a number of
assumptions that might be difficult to verify for media of high
contrast and high frequency permeability.  Nevertheless, the results in 
\cite{Willems_2013,Willems_2014}, in our opinion, present the state of the art of
preconditioning of such problems.

The numerical simulation of processes in media of high frequency and
high contrast represents a great challenge since it leads to
ill-conditioned symmetric but indefinite system of linear
equations. In our opinion, its efficient preconditioning is not fully
mastered.  Here we present a step in this direction for Brinkman/Darcy
models.  The main objectives and contributions of our paper are:
\begin{enumerate}
\item To discuss and present a unified solution methodology 
for computer simulation of flows in porous media modeled by
  Darcy and Brinkman equations. 
  Using this methodology, one may set up natural experiments with
  highly heterogeneous media in order to compare and analyze the
  numerical simulations in the framework of a mathematical modeling
  tool.
\item To show the efficiency of the developed preconditioner for
  solving very large systems of linear equations arising from the
  finite element approximation of the Darcy and Brinkman equations and
  to demonstrate via various tests the robustness of the method with
  respect to both, the mesh step-size and the high contrast high
  frequency porous media.
\item To experiment with various two- and three-dimensional synthetic
  test problems that are used by flow in porous media community and
  show the capabilities of numerical simulation methodology for
  solving relevant applied problems.
\end{enumerate}

The remainder of this paper is organized as follows: In Section
\ref{sec:Brinkman}, we provide a detailed description of the problems
under consideration as well as the necessary
notation. 
Further, we provide a description of a DG discretization of Brinkman's
equations. In Section~\ref{sec:multigrid} we outline the derivation of
the multigrid algorithm for Brinkman's and Darcy's equations.
In Section \ref{sec:numerics} we present numerical experiments with
benchmark problems that demonstrate the capabilities and the
robustness of our method. In the final Section \ref{sec:conclusions} 
we present the main conclusions of the paper.

\section{Problem Formulation and Notation}\label{sec:Brinkman}

\subsection{Notation of spaces}
We use the standard notation for spaces of scalar and vector 
functions defined on a bounded domain $\qOmega\subset\qbR^d$ ($d=2,3$)
with polyhedral boundary having the outward unit normal vector $\qbn$.
$L^2_0(\qOmega)\subset L^2(\qOmega)$ denotes the space of square
integrable functions with mean value zero and $H^1(\qOmega)^\qn$,
$H^1_0(\qOmega)^\qn$, and $L^2(\qOmega)^\qn$ denote the spaces of
vector-valued functions with components in $H^1(\qOmega)$,
$H^1_0(\qOmega)$, and $L^2(\qOmega)$, respectively. Furthermore, we set
\begin{alignat*}2
  H^{\text{div}}( \qOmega) &=\{\qbv \in L^2(\qOmega)^\qn
  &\;:\;& \div \qbv \in L^2(\qOmega)\}.
\end{alignat*}
The spaces 
$H^{\text{div}}( \qOmega) $  and $H^1(\qOmega)^d$ are 
equipped with the following inner products
$$
(\qbu,\qbv)_{H^{\text{div}}( \qOmega)}=\dint{\qOmega}{} ( \div \qbu \, \div \qbv + \qbu \cdot \qbv) \qdbx, 
\quad \mbox{and} \quad 
(\qbu,\qbv)_{H^{1}( \qOmega)}=\dint{\qOmega}{} ( \grad \qbu  : \grad \qbv + \qbu \cdot \qbv) \qdbx,
$$
 and corresponding norms
$
\dnorm{\qbv}{H^{\text{div}}(\qOmega)}^2 =(\qbu,\qbv)_{H^{\text{div}}(\Omega)}\ \  \mbox{and} \ \  
\dnorm{\qbv}{H^1(\qOmega)}^2 =(\qbu,\qbv)_{H^{1}( \qOmega)}.
$
%

\subsection{Problem formulation}
\def\bu{\boldsymbol{u}}
We consider the Brinkman/Darcy model 
for the macroscopic pressure $p$
and the fluid velocity $\qbu=(\qu_1,\ldots,\qu_\qn)$: 
\begin{equation}
  \label{eq:Brinkman}
-\qmu\laplace\qbu+\qkappa\qbu +\grad\qp= \qbf, \quad   \div\qbu = 0     \quad \text{on }\qdOmega.
\end{equation}
Note that formally one gets the Darcy model by setting $\qmu=0$.  To
this system we add proper boundary conditions. We shall consider the
simplest ones
\begin{equation}\label{BC}
\qbu = \qbg \quad \left ( \text{or }\qbu \cdot \qbn = \qbg\cdot \qbn
  \text{ if } \mu=0 \right ) \quad 
\text{on }\qdOmega.
\end{equation}

We assume that the boundary data $\qbg\in H^{\frac12}(\qdOmega)^\qn$
(or $\qbg\cdot \qbn \in H^{-\frac12}(\qdOmega)$ for the Darcy model)
and satisfy the compatibility condition
$$
\int_{\qdOmega} \qbg \cdot \qbn \qds =0.
$$

Problem \eqref{eq:Brinkman} has unique  solutions 
$(\qbu,\qp)$ in $H^1(\qOmega)^\qn\times L^2_0(\qOmega)$ or,
if $\qmu=0$, in $H^{\text{div}}(\qOmega)\times L^2_0(\qOmega)$.
The smoothness of the solutions of these problems can be studied by the methods 
developed e.g. in \cite{Grisvard}. However, due to the media character, high contrast and high frequency,
this problem has limited solution regularity, which depends on the coefficients jumps and their arrangement
in rather unfavorable way. 

\subsection{Discontinuous Galerkin (DG) Finite Element Method}
For the finite element approximation of \eqref{eq:Brinkman} we shall
employ Raviart-Thomas finite elements. To this end, we partition the
domain $\Omega$ into rectangular and hexahedral cells of size $h$  in two
and three dimensions, respectively, and denote this partitioning by
$\qTcalH $.
Further, we shall need a notation for the set $ \qiEcalH$ of the
internal and $\qdEcalH$ for the boundary interfaces/edges of the
finite element partition.

Now, let $\qWscrh \subset L^2_0(\qOmega)$ be the space of discontinuous 
functions that are polynomial of degree $k>0$ in each variable $x_i$, $i=1, \dots, d$ and
$\mathscr{V}_h \subset H^{\text{div}}(\qOmega)$ 
be a finite element
space based on Raviart-Thomas ($RT_k$) elements of degree $k \ge 0$ 
(cf.\ e.g.\ \cite[pages 120--130]{BF91}).
%
Since 
$\qVscrH \not\subseteq H^1(\Omega)^d$ 
the gradient of functions in $ \qVscrH$ will be defined only in the
interior of each element.  Nevertheless, the theory of DG FEM developed in 
~\cite{CockburnKanschatSchoetzau05} for $k \ge 1$ shows that the scheme is stable and 
has optimal approximation properties.
In the special case of lowest order Raviart-Thomas
elements ($RT_0$ or $k=0$), the gradient is not approximated consistently
and the theory in~\cite{CockburnKanschatSchoetzau05} fails. Nevertheless, it was
shown in~\cite{Kanschat08mac}, that in this case we get the MAC scheme of
Harlow and Welch~\cite{HarlowWelch65}, for which 
an alternative theory exists, e.g. \cite{Nicolaides_MAC}.

If $T_1$ and $T_2$ 
are two mesh cells with a common face $F$ and  $\qbv_1$ and $\qbv_2$ are the traces of a 
function $\qbv \in \qVscrH$ on $F$  from $T_1$ and $T_2$, respectively, then
we define the average and the jump by
\begin{gather}
  \mvl{ \qbv } = \frac12(\qbv_1 + \qbv_2) \quad \mbox{and} \quad 
\jump{\qbv} = \qbv_1 - \qbv_2  \quad \mbox{for} \quad x \in F, \, \,   F  \in \qiEcalH.
\end{gather}
For the jump to be defined consistently, we assume that there is a global ordering of the finite elements 
so that $T_1$ is the element with a smaller number and $T_2$ is the one with higher number.
It is useful to have the same notation for faces that are on the boundary $\partial \Omega$,
namely, $ \mvl{ \qbv } = \qbv$ and $ \jump{\qbv} = \qbv$ on $ F \subset \partial \Omega$.

Now we introduce the 
following notations for integrals of vector functions  $\qbu_h, \qbv \in  \qVscrH$
over the domain $\Omega$ and over the finite element
interfaces:
$$
(\nabla \qbu_h, \nabla \qbv)_{\qTcalH} 
= \sum_{T \in \qTcalH} \int_T \nabla \qbu_h : \nabla \qbv dx,  \quad  
\langle \jump{\qbu_h}, \jump{\qbv} \rangle_{\qiEcal} 
= \sum_{F \in \qiEcal} \int_F  \jump{\qbu_h} \cdot \jump{\qbv} \, ds
$$
and similarly for $\langle \mvl{\qbu_h}, \jump{\qbv} \rangle_{\qiEcal} $. 
Further we will use also the notation
$$
\forme(\qbu_h,  \qbv)_{\qiEcalH}  
:=\sum_{F \in \qiEcalH} \int_F  \qbu_h \cdot \qbv \, ds 
\quad \mbox{and} \quad
\forme(\qbu_h,  \qbv)_{\qdEcalH}  
:=\sum_{F \in \qdEcalH} \int_F  \qbu_h \cdot \qbv \, ds = \int_{\partial \Omega} \qbu_h \cdot \qbv \, ds.
$$
Thus, for function $ \qbu_h, \qbv \in V_h$
we have the following obvious relations:
$$ \forme( \jump{\qbu_h}, \jump{\qbv})_{\qiEcal}   
= \forme( \jump{\qbu_h}, \jump{\qbv})_{\qiEcalH} + \forme(\qbu_h, \qbv)_{\qdEcalH}, \quad 
\forme( \mvl{ \qbu_h},\jump{\qbv} )_{\qiEcal} =
 \forme( \mvl{ \qbu_h},\jump{\qbv} )_{\qiEcalH} 
   + \forme(\qbu_h, \qbv)_{\qdEcalH}.
$$

We use a discontinuous Galerkin (DG) finite element method to
discretize problem \eqref{eq:Brinkman}.
Following~\cite{CockburnKanschatSchoetzau07,KanschatRiviere10}, we
define the interior penalty bilinear form for $\qbu_h , \, \qbv \in \qVscrH$ 
and for $\sigma_h =\sigma/h$, $ \sigma >0$ a sufficiently large stabilization parameter 
\begin{gather}
  \begin{split}
    \label{eq:ipbrinkman}
    a_h(\qbu_h,\qbv) =& \mu (\nabla \qbu_h, \nabla \qbv)_{\qTcalH} + (\kappa \qbu_h, \qbv)_{\qTcalH}
    + 2\mu{\sigma_h}\forme( \jump{\qbu_h}, \jump{\qbv})_{\qiEcal}
    \\
    &- \mu \forme( \mvl{\nabla \qbu_h}\cdot \qbn,\jump{\qbv} )_{\qiEcal} 
   - \mu \forme(\mvl{\nabla \qbv }\cdot \qbn, \jump{\qbu_h})_{\qiEcal},
  \end{split}
\end{gather}
where $ \mvl{\nabla \qbv}$ is the average of the $d \times d$ matrix 
$ \nabla \qbv$  and
$\qbn$ is a unit normal vector, fixed for each internal face $F \in \qiEcalH$  and 
pointing outward $\Omega$
for a boundary face $ F \subset \partial \Omega$. 

The discrete weak formulation of \eqref{eq:Brinkman} reads now: find $(\qbu_h,p_h)
\in \qVscrH \times \qWscrh$, such that 
\begin{gather}
  \label{eq:Brinkman-weak}
\A_h(\qbu_h, p_h; \qbv, q):=  a_h(\qbu_h,\qbv)
  - (p_h, \div \qbv)
  - \form(q, \div \qbu_h)
  = \mathcal{F}(\qbv, q) , \quad \forall (\qbv, q)  \in \qVscrH \times  \qWscrh,
\end{gather}
where $ \mathcal{F}(\qbv, q)$
contains the right-hand side $\qbf$ and the Dirichlet boundary data $\qbg$:
\begin{equation}\label{eq:BCh}
 \mathcal{F}(\qbv, q) =(\qbf, \qbv) + 2\mu \sigma_h\forme(\qbg, \qbv)_{\qdEcalH} - \mu  \forme(\nabla\qbv\cdot \qbn, \qbg)_{\qdEcalH}.
\end{equation}
The discussion on the existence and uniqueness of  solutions for the Brinkman 
system ~\eqref{eq:Brinkman-weak}, \eqref{eq:BCh} 
can be found in ~\cite{IlievLazarovWillems11,MaoDissertation}. 
%
We get Darcy case by taking $\mu =0$
and 
$ \mathcal{F}(\qbv, q) =(\qbf, \qbv) + \forme(\qbg \cdot \qbn, q)_{\qdEcalH} $.



\section{The multigrid method}
\label{sec:multigrid}

The preconditioner is based on the geometric multigrid method that uses a hierarchy of meshes 
$\mathcal{T}_{j}$, $j=0, \dots,J$. The mesh $\mathcal{T}_{0}$ is the initial coarse grid 
partition and $\mathcal{T}_{J} \equiv \mathcal{T}_{h}$ is the finest grid, where the approximation~\eqref{eq:Brinkman-weak}--\eqref{eq:BCh} is set. The mesh hierarchy is defined recursively,
such that the cells of $\mathcal{T}_{j+1}$ are obtained by splitting
each cell of $\mathcal{T}_j$ into $2^d$ children by connecting the
edge, face, and cell midpoints (refinement).  These meshes are nested
in the sense that every cell of $\mathcal{T}_j$ is equal to the union
of its four children (eight in 3-D case). 
We define the mesh size $h_j$ as the maximum of the diameters of the
cells of $\mathcal{T}_j$. 

\subsection{Nested FE spaces}


Applying the finite element spaces introduced above to the meshes $\qTcalj$, we obtain a nested sequence of spaces
\begin{gather*}
  \arraycolsep5pt
  \begin{array}{rcccccl}
    \mathscr{V}_0 &\subset& \mathscr{V}_1 &\subset& \dots &\subset& \mathscr{V}_J,
    \\
    \mathscr{Q}_0 &\subset&  \mathscr{Q}_1 &\subset& \dots &\subset&  \mathscr{Q}_J.
    \\
     \mathscr{V}_0 \times  \mathscr{Q}_0 =: \mathscr{X}_0  &\subset&  \mathscr{X}_1 &\subset& \dots
    &\subset& \mathscr{X}_J := \mathscr{V}_J \times  \mathscr{Q}_J. \\
  \end{array}
\end{gather*}

The nestedness of the spaces implies that there is a sequence of
natural injections $\I_j: \mathscr{X}_j \to \mathscr{X}_{j+1}$ of the form
$\I_j(v_j, q_j) = (I_{j,u}v_j, I_{j,p}q_j)$, such
that
$$ 
  I_{j,u}: 
\mathscr{V}_j \to \mathscr{V}_{j+1},
  \qquad
  I_{j,p}:\;  \mathscr{Q}_j \to \mathscr{Q}_{j+1}.
$$ 
From the cochain complex properties of the Raviart-Thomas spaces, we
have additionally, that the divergence free subspaces
$\mathscr{V}_j^0$ are nested and thus there holds
$I_{j,u}: \mathscr{V}_j^0 \to \mathscr{V}_{j+1}^0$.  The
$L^2$-projection from $\mathscr{X}_{j+1} \to \mathscr{X}_j$ is defined
by $\I^t_j(v_j, q_j) = (I^t_{j,u}v_j,I^t_{j,p}q_j)$ with
\begin{equation}
  \label{eq:3}
  \begin{split}
  \form(v_{j+1}-I^t_{j,u}v_{j+1}, w_j) &= 0 \quad \forall w_j
  \in \mathscr{V}_j 
\\
  \form(q_{j+1}-I^t_{j,p}q_{j+1}, r_j) &= 0 \quad \forall r_j
  \in \mathscr{Q}_j.    
  \end{split}
\end{equation}


\subsection{The variable V-cycle algorithm}
\label{sec:V-cycle}


Due to the multilevel structure of the spaces $\qVscrj \times  \qWscrj$, $j=0, \dots,J$,
on each grid level $j$ 
the weak formulation to find $(\qbu_j, p_j) \in \qVscrj \times  \qWscrj$ such that
\begin{gather}\label{ML-problem}
\A_j(\qbu_j, p_j; \qbv, q)
= \form(\qbf,\qbv) \,\, \,
\forall (\qbv, q)  \in \qVscrj \times  \qWscrj,
\end{gather}
is rewritten in an algebraic form $\AAA_j x = b_j$, where $x_j =(u_j, p_j)$ and
$b_j = (f_j,g_j)$. Then the discretization ~\eqref{eq:Brinkman-weak},
\eqref{eq:BCh} on the finest grid has the following algebraic form
\begin{gather*}
  \AAA_J x_J = b_J.
\end{gather*}
The system has saddle point structure with respect to the variables 
for velocity and pressure. We solve it by a preconditioned GMRES iteration.
Due to the elliptic structure of the Brinkman operator, we choose the 
multilevel preconditioner $\B_J$ defined in the remainder of this section.
We assume that the coarse mesh $\mathcal{T}_{0}$ has small number of degrees of freedom
so we can afford direct solution of the system $ \AAA_0 x_0 = b_0$.

We define the preconditioner
$\B_j : \mathscr{X}_j \rightarrow \mathscr{X}_j $ recursively. Let
$\R_j$ being a suitable smoother, as described below and let
$m(j)$ be the number of smoothing steps on level $j$.  Let
$\B_0 =\AAA_0^{-1}$.  For $j= 1,\ldots,J$ define the action of $\B_j$ on
a vector $b_j:=(f_j, g_j) \in \mathscr{X}_{j} $ as follows:
\begin{subequations}
\begin{enumerate}
\item Pre-smoothing: begin with $x_0 =  0 \in \mathscr{X}_{j}$ 
and compute $x_i \in  \mathscr{X}_{j}$ with $i= 1, \dots, m(j)$ by
\begin{gather*}
  x_i  = x_{i-1} + \R_{j}
  \left( b_j - \AAA_j  x_{i-1},
  \right).
\end{gather*}
  
\item Coarse grid correction:
\begin{gather*}
x_{m(j)+1}    =     \B_{j-1} \I_{j-1}^t \left(b_j - \AAA_j \, x_{m(j)}    \right).
\end{gather*}

\item Post-smoothing: for $ i= m(j)+2, \dots, 2m(j)+1$ compute
\begin{gather*}
 x_i    = x_{i-1}    + \R_{j}  \left(b_j - \AAA_j \, x_{i-1},    \right).
\end{gather*}
\item Assign: 
$ 
    \B_j  b_j = x_{2m(j)+1}.
$ 
\end{enumerate}  
\end{subequations}

We distinguish between the standard V-cycle with $ m(j) =m(J)$ and the
variable V-cycle with $m(j) =m(J) 2^{J-j} $, where the number $m(J)$
of smoothing steps on the finest level is a parameter. In our
numerical examples below we use the variable V-cycle with $m(J)= 2$.
We refer to $\B_J$ as the V-cycle preconditioner of $\AAA_J$,
independent of the choice of $m(j)$.

\subsection{Overlapping Schwarz smoothers}
\label{sec:smoother}

\begin{figure}[tp]
  \centering
                      \includegraphics[height = 3cm, width=3cm]{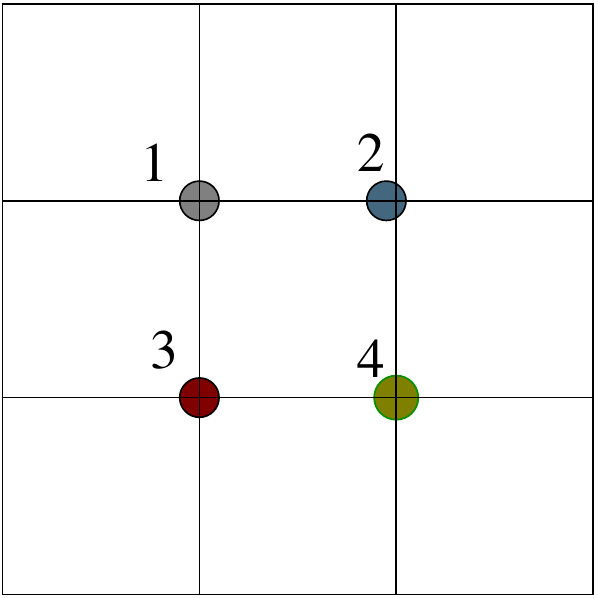}\hspace{1cm}
                      \includegraphics[height = 3cm, width=3cm]{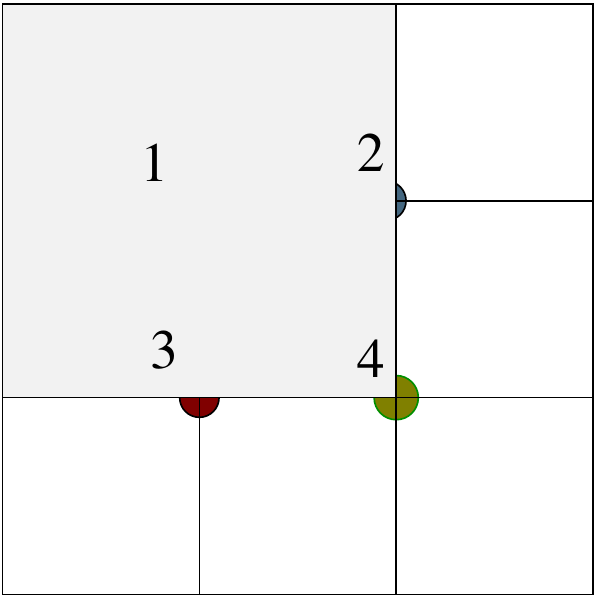}\hspace{1cm}
                      \includegraphics[height = 3cm, width=3cm]{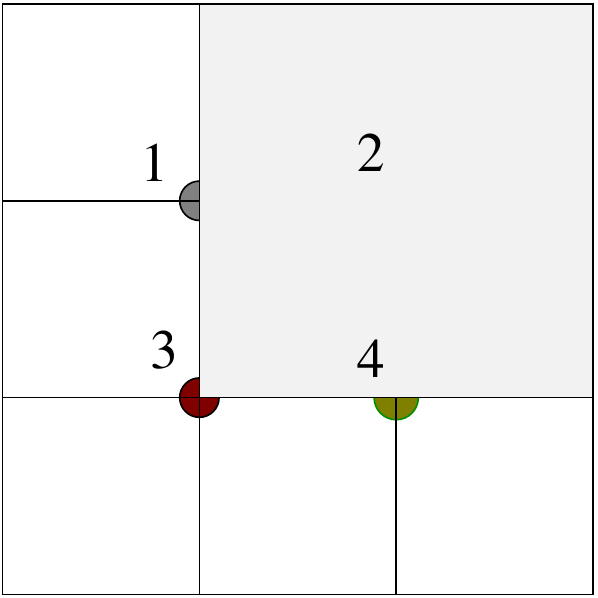}
                      \caption{Part of the domain consisting of patches around four vertices (left). In the center we have
                                     a patch around the vertex 1 and on the right -- a patch around the vertex 2}\label{patch}
\end{figure}

In this subsection, we define a class of smoothing operators $\R_j$ based 
on a subspace decomposition of the space  $ \mathscr{X}_j$.
Let
$\mathcal{N}_j$ be the set of vertices in the triangulation $\mathcal{T}_j$, and let
$\mathcal{T}_{j,\upsilon}$ be the set of cells in $\mathcal{T}_j$ sharing the
vertex $\upsilon$ (see Figure~\ref{patch}). They form an overlapping covering 
with $N_j$ patches ($N_j>0$),
denoted by $\{\Omega_{j,\upsilon}\}_{\upsilon = 1}^{N_j}$.

The subspace  $ \mathscr{X}_{j,\upsilon} = \mathscr V_{j,\upsilon} \times
\mathscr Q_{j,\upsilon}$
 consists of the functions in $ \mathscr{X}_j$ with support
in $\Omega_{j,\upsilon}$. Note that this implies homogeneous slip
boundary conditions on $\partial\Omega_{j,\upsilon}$ for the
velocity subspace $\mathscr V_{j,\upsilon}$ and zero mean value on
$\Omega_{j,\upsilon}$ for the pressure subspace
$\mathscr Q_{j,\upsilon}$. The Ritz projection $\mathcal P_{j,\upsilon}:
	 \mathscr{X}_j \to  \mathscr{X}_{j,\upsilon}$ is defined by the equation
$  \A_j\form(\mathcal P_{j,\upsilon} x_{j}, y_{j,\upsilon})
  = \A_j\form(x_{j}, y_{j,\upsilon}),
  \quad \forall  y_{j, \upsilon} \in  \mathscr{X}_{j,\upsilon}.
$ 
This equation 
involves solving a local saddle point problem on the patch $\Omega_{j,\upsilon}$.
Note that each cell belongs to not more than four (eight in case of 3-dimentions)
patches $\Omega_{j,\upsilon}$, one for each of its vertices.


Following~\cite{TW05} we define the symmetric multiplicative Schwarz smoother $\R_{j}$,
associated with the spaces $ \mathscr{X}_{j,\upsilon}$,   by
\begin{equation}
\label{eq:smoother-Multiplicative}
  \R_{j} = (\mathcal I - \mathcal E_j\mathcal E^{\ast}_j) \AAA^{-1}_j
  \quad \mbox{with} \quad 
  \mathcal E_j =
  \left(\mathcal I - \mathcal P_{j,N_j}\right)
  \cdots
  \left(\mathcal I - \mathcal P_{j,1}\right),
\end{equation}
where  $\mathcal{E}^{\ast}_j$ is the $\AAA_j$-adjoint of $\mathcal E_j$.
Note that $\AAA^{-1}_j$ is never computed, it is part only of the theoretical justification and 
formulation of the preconditioner, see, for more details \cite{TW05}.

\section{Numerical Experiments}\label{sec:numerics}

We perform numerical tests on several benchmark configurations. First,
we introduce several representative permeability fields.
Further, these are complemented by permeability 
fields from the SPE10 benchmark~\cite{ChristieBlunt01}.

\subsection{Setting the parameters in the computational procedure.}
In all numerical experiments, we use the variable
V-cycle with two smoothing steps on the finest level and the
symmetric, multiplicative Schwarz smoother $\R_{j}$.  The
penalty parameter $\sigma_h$ in \eqref{eq:ipbrinkman}, which is denoted by $\sigma_j$
in the multilevel setting \eqref{ML-problem}, 
is chosen as
$\sigma_j=(k+1)(k+2)/h_j$, where $h_j$ is the mesh-size on level   $j$. 
%
The stopping criterion in the GMRES solver is set to reduce the
Euclidean norm of the initial residual by $10^{-6}$.  

For setting up the discrete systems
$ 
  \A_j\left( x_j, y_j \right )  = \mathcal{F}(y_j), \, \,  x_j=(\qbu_j,p_j), \,\, \,  y_j=(\qbv_j,q_j)
$ 
on each level $j$ we also need 
the permeability on each grid level. In our implementation we use
arithmetic averaging of permeabilities. Other choices are also possible,
see, e.g.~\cite{MaierThesis}.

\subsection{Tests on some characteristic permeability fields}
Our first numerical tests are done 
on a number of  two dimensional permeability fields generated on a
$128 \times 128$ fine mesh as shown in Figure~\ref{fig:set1}.  Here,
we have used the following convention: the red color denotes regions
of media with high permeability, while the blue color
corresponds to media with low  permeability (in black and white  -- 
gray represents high permeability, while dark means low permeability).  
On Figure~\ref{fig:set1} in the red regions $\kappa(x)=1$ 
and in the blue regions $\kappa(x)=1/{\boldsymbol \kappa}$.
These
have been used as representative test examples by other
authors, see, e.g.  \cite{LefebvreBanhartDunand08,
  PopovBiEfendievEwingQinLi07}.

Figure~\ref{fig:set1}(A) shows the permeability of a two-dimensional 
cut of open foams of porous media (similar to those shown in Figure~\ref{fig:introduction}(A)). 
In this example, low  permeability regions are mostly disconnected. 
Figures \ref{fig:set1}(B) and (C) 
present porous media with isolated and globally connected  inclusions  
with low   permeability. 
Figure~\ref{fig:set1}(C) represents a media with some long connected  inclusions that 
mimic bonded fibrous porous materials.

\subsubsection{Darcy vs.~Brinkman mathematical models}

\begin{figure}[tp]
        \centering
        \begin{subfigure}[H]{0.32\textwidth}
                    \includegraphics[width=\textwidth,height=\textwidth]{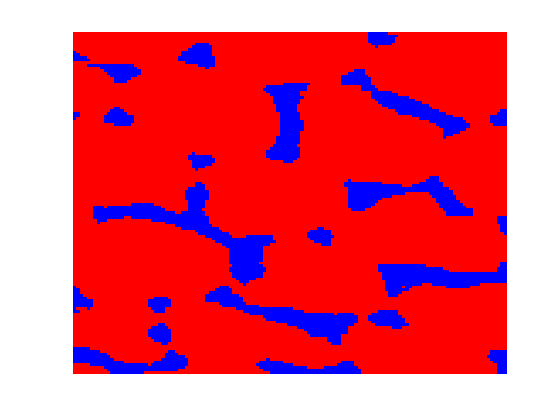}
                \caption{Open foam}   
        \end{subfigure}
        \hfill
        \begin{subfigure}[H]{0.32\textwidth}
                    \includegraphics[width=\textwidth,height=\textwidth]{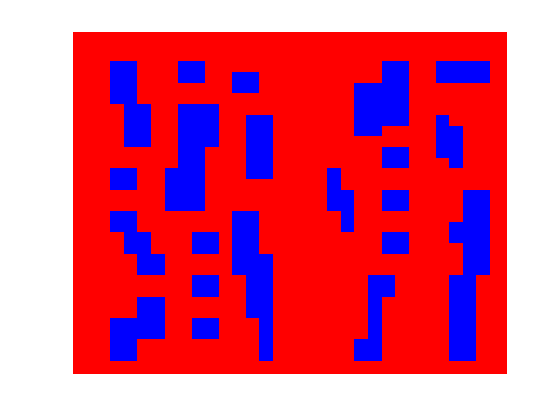}
                \caption{{Inclusions} } 
        \end{subfigure}
        \hfill
        \begin{subfigure}[H]{0.32\textwidth}
                    \includegraphics[width=\textwidth,height=\textwidth]{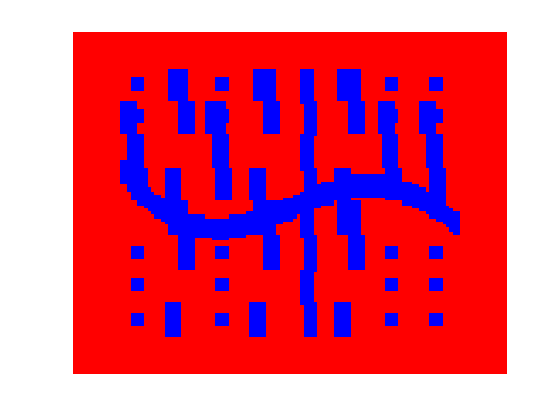}
                \caption{{Connected inclusions} }   
        \end{subfigure}
      \caption{Three distributions of permeability: red regions denote 
        high permeability, blue regions mark low permeable inclusions}
      \label{fig:set1}
\end{figure}

In order to examine the effects of numerical modeling on fluid flow in 
porous media, we have run tests using Darcy and Brinkman models on three 
different sets of permeability distribution shown in
Figure~\ref{fig:set1}. The flow is driven by boundary values $\qbg =
(1,0)$ and $\qbf=\q0$.
For Brinkman's model we use $\mu = 10^{-2}$ and, obviously, $\mu=0$
for Darcy's model. The contrast $\boldsymbol{\kappa}$ is in the range $10^4$  --
$10^6$. 

  \begin{figure}[tp]
        \centering
        \begin{subfigure}[H]{0.3\textwidth}
          \includegraphics[width=\textwidth]{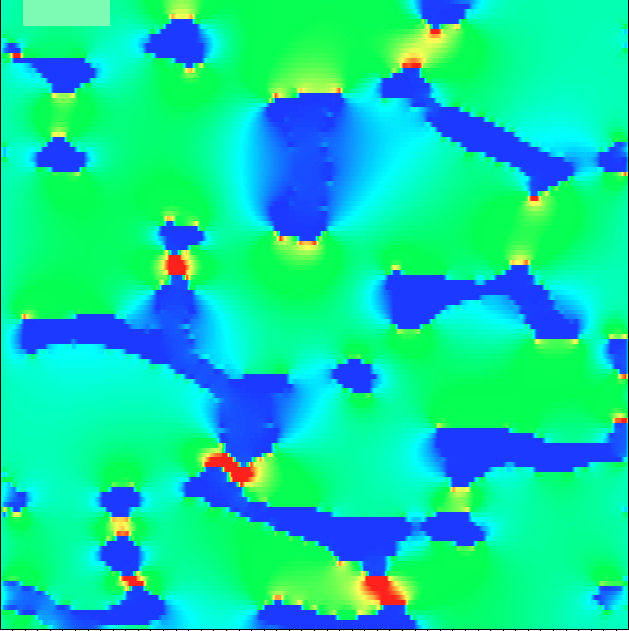}
          \caption{Open foam, Darcy}
        \end{subfigure}
        \hfill
        \begin{subfigure}[H]{0.3\textwidth}
        \includegraphics[width=\textwidth]{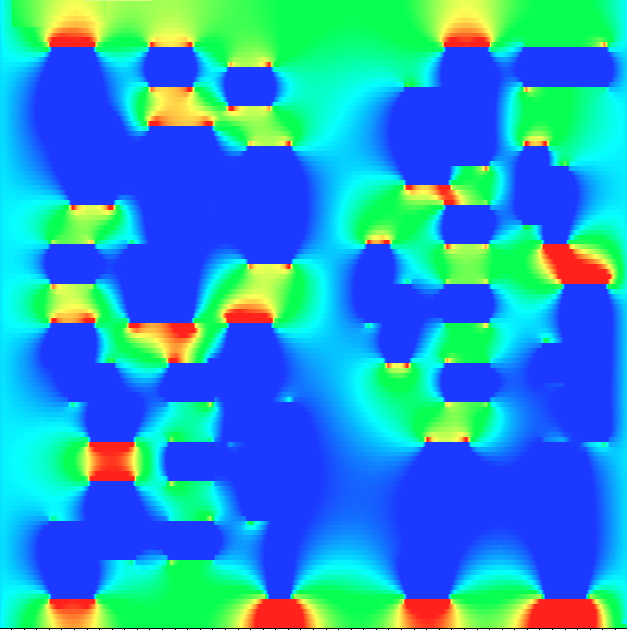}
          \caption{Inclusions, Darcy}
        \end{subfigure}
        \hfill
        \begin{subfigure}[H]{0.3\textwidth}
        \includegraphics[width=\textwidth]{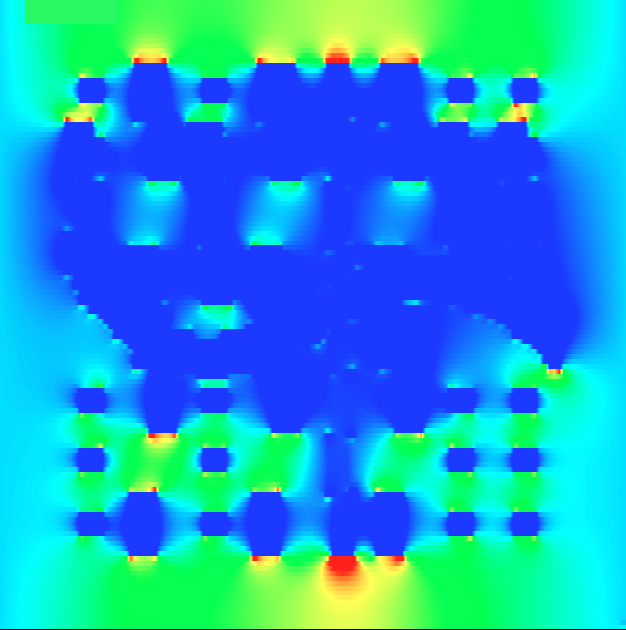}
          \caption{Con. inclusions, Darcy}
        \end{subfigure}
        \\[2mm]
        \begin{subfigure}[H]{0.3\textwidth}
        \includegraphics[width=\textwidth]{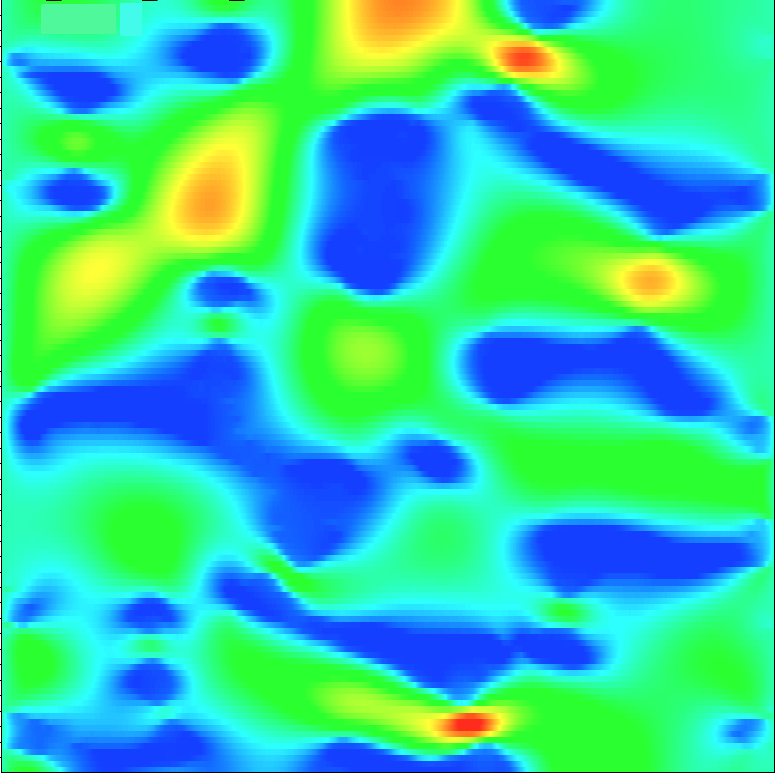}
          \caption{Open foam, Brinkman}
        \end{subfigure}
        \hfill
        \begin{subfigure}[H]{0.3\textwidth}
        \includegraphics[width=\textwidth]{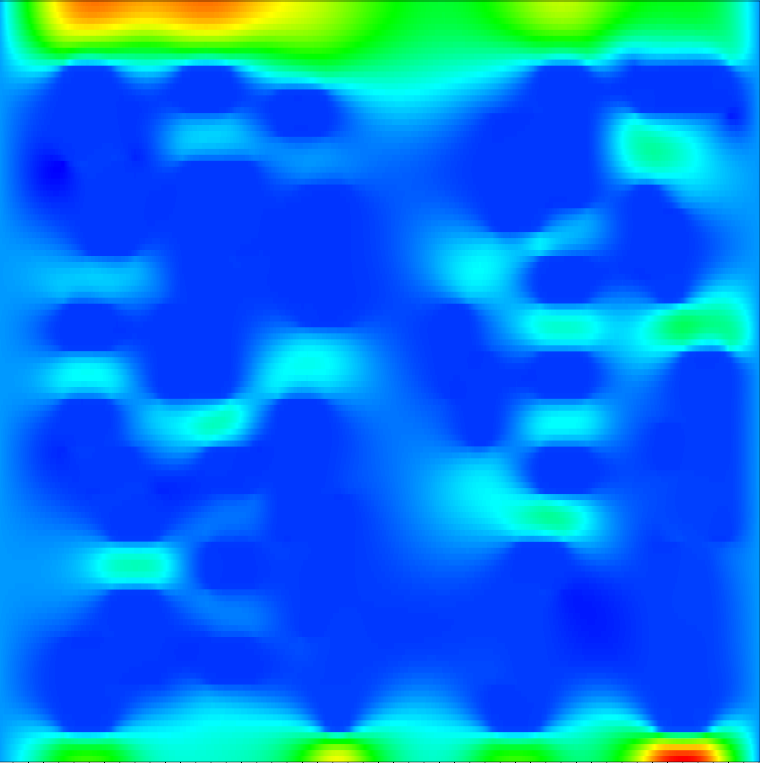}
          \caption{Inclusions, Brinkman}
        \end{subfigure}
        \hfill
        \begin{subfigure}[H]{0.3\textwidth}
        \includegraphics[width=\textwidth]{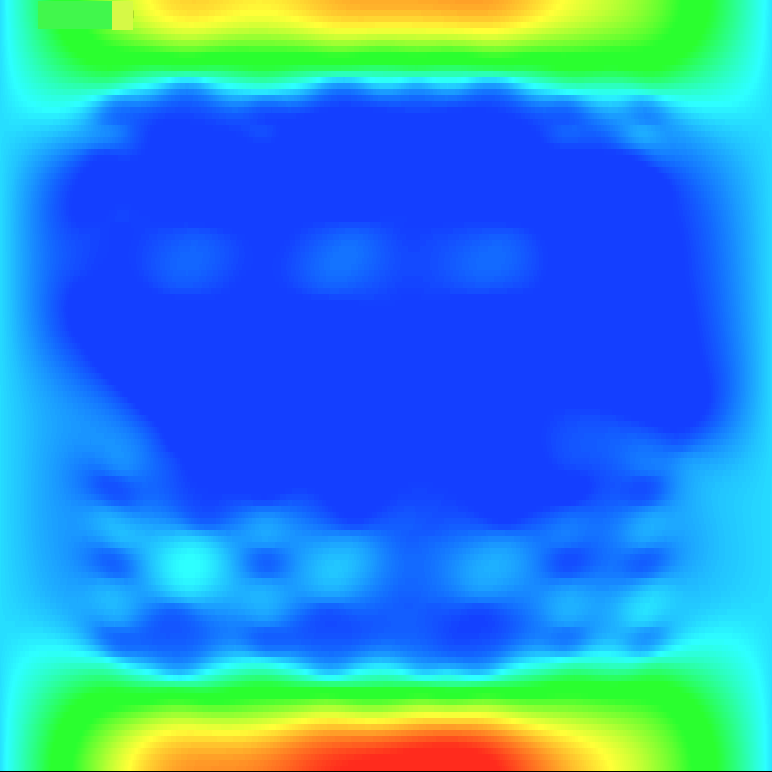}
          \caption{Con. inclusions, Brinkman}
        \end{subfigure}

        \caption{
         Velocity component $u_1$  for permeability 
        distributions of Figure~\ref{fig:set1} for Darcy and Brinkman equations with contrast 
        $\boldsymbol{\kappa}=10^6$: for open foam, the maximal value (red) is 3.47, minimal value (blue) is negative 0.13; for inclusions, the maximal value (red) is 3.89, minimal value (blue) negative 0.26; for connected inclusions, the maximal value (red) is 4.37, the minimal value (blue) negative 0.09. }
        \label{fig:comp1}
\end{figure}  

The computational results for the $u_1$ component of  the velocity
are presented in Figure~\ref{fig:comp1}: (A) -- (C) are for Darcy and (D) -- (F) --
for Brinkman models.  The velocity is presented in the same color scale so one 
can observe similarity in the patterns of the flow. 
%
One can also see
different flow behaviors such as substantial variation of velocity magnitude.

One can make several immediate observations: 
(1) there are apparent diffusion effect in the flow governed by Brinkman model; 
(2) the connected inclusions, Figure~\ref{fig:set1}(C), are
blocking the Brinkman flow in a more profound way compared with the Darcy case;
(3) the most visible differences are in the upper horizontal boundary,
where no-vertical flow for both Brinkman and Darcy are imposed; in the case of Brinkman model,
most of the fluid is flowing through the region adjacent to the horizontal boundaries
since in the rest the flow  is blocked  by the inclusion that connects other long inclusions.

In addition, one visually observes a substantial difference between
solutions produced by Darcy and Brinkman models while changing the
permeability contrast. In order to quantify these visual difference,
we computed the $H^1$-norm of these solutions. 

  More numerical experiments on these sets of data and some other permeability fields can be found in 
  \cite{MaoDissertation}. Based on these numerical experiments we observed that
for Brinkman model (with permeability shown in Figure~\ref{fig:comp1} (B), disconnected
inclusions), the difference between the solutions
obtained for $\boldsymbol{\kappa}=10^4$ and $\boldsymbol{\kappa}=10^5$
is $22\%$, while for $\boldsymbol{\kappa}=10^5$ and $\boldsymbol{\kappa}= 10^6$ it
is only $8\%$. Further increase of the contrast does not lead to substantial change of the solution. 
For Darcy and
Brinkman model with the same permeability contrasts shown in Figure~\ref{fig:comp1}, we found that the difference between the solutions is up to $600\%$. 

Based on our numerical experiments we conclude that the two numerical
models produce results which are robust with respect to the
parameters.  The choice of the model to be used in a particular
practical problem will depend on physics involved and the agreement
with certain natural experiments.

\subsubsection{Performance of the preconditioner on two-dimensional Darcy and Brinkman models}
\label{darcy solver}

We also test the robustness (with respect to the contrast $\boldsymbol{\kappa}$
and the mesh-size $h$) of the preconditioner
for Darcy and Brinkman models for flows in heterogeneously porous media 
on three sets of heterogeneous permeabilities shown in Figure~\ref{fig:set1}
in two cases:
\begin{description}
\item[case (a)] 
 $\kappa(x) = 1$  in the blue regions (of high ~permeability, channels-like formations) and 
$\kappa(x)=10^{-4}, 10^{-5}, 10^{-6}$ in the red regions;
\item[case (b)] 
we reverse the roles of low and high permeability by taking $\kappa(x) = 1$  in the red region (of low ~permeability) and 
$\kappa(x)=10^{-4}, 10^{-5}, 10^{-6}$ in the blue regions.
\end{description}
  \begin{table}[h]
  \center
  \begin{tabular}{ | c | c c c | c c c | c c c |}
    \hline
    & \multicolumn{3}{c|}{Open foam} 
    & \multicolumn{3}{c|}{Inclusions}
    & \multicolumn{3}{c|}{Con. inclusions}    
    \\\hline
    \backslashbox{$h^{-1}$}{$\kappa$}
    &$10^4$&$10^5$&$10^6$
    &$10^4$&$10^5$&$10^6$
    &$10^4$&$10^5$&$10^6$
    \\ \hline
    &\multicolumn{9}{c|}{Darcy}
      \\\hline
    128&17&19&20&40&51&60&44&61&73
    \\   \hline
    256&16&17&19&35&45&54&41&55&65
    \\ \hline
    512&15&16&18&33&42&50&38&51&59
    \\ \hline
    & \multicolumn{9}{c|}{Brinkman}
    \\\hline
    128&18&27&36&16&22&38&16&23&36
    \\   \hline
    256&16&24&32&15&20&32&14&20&31
    \\ \hline
    512&15&22&31&14&19&29&13&19&29
    \\\hline
  \end{tabular}
  \caption{Number of GMRES iterations needed to reduce the residual
    by $10^{-6}$ for Darcy and Brinkman models for permeability
    shown in Figure~\ref{fig:set1}} 
  \label{geom_set_1}
\end{table}
  \begin{table}[h]
  \center
  \begin{tabular}{ | c | c c c | c c c | c c c |}
    \hline
    & \multicolumn{3}{c|}{Open foam} 
    & \multicolumn{3}{c|}{Inclusions}
    & \multicolumn{3}{c|}{Con. inclusions}    
    \\\hline
    \backslashbox{$h^{-1}$}{$\kappa$}
    &$10^4$&$10^5$&$10^6$
    &$10^4$&$10^5$&$10^6$
    &$10^4$&$10^5$&$10^6$
    \\ \hline
    &\multicolumn{9}{c|}{Darcy}
      \\\hline
    128&55&67&78&76&81&88&68&74&81
    \\   \hline
    256&48&61&71&68&75&82&63&69&76
    \\ \hline
    512&45&58&67&64&70&78&61&64&70
    \\ \hline
    & \multicolumn{9}{c|}{Brinkman}
    \\\hline
    128&48&53&58&65&73&79&68&73&77
    \\   \hline
    256&43&49&54&61&68&72&63&69&72
    \\ \hline
    512&39&46&49&57&63&68&59&63&69
    \\\hline
  \end{tabular}
  \caption{Number of GMRES iterations needed to reduce the residual
    by $10^{-6}$ for Darcy and Brinkman models for 
    permeability shown in Figure~\ref{fig:set1}, with low and high permeability reversed} 
  \label{geom_set_reverse}
\end{table}
In each test of Tables \ref{geom_set_1} and \ref{geom_set_reverse}  for a fixed 
permeability contrast (in a column), we observe 
iteration counts practically independent of the mesh size $h$. 
This shows that the preconditioner is uniform (possibly optimal) with respect to the
mesh-size.

For a fixed fine-grid (in a row), we report the iteration count while
changing the permeability in both, Darcy and Brinkman models.
For a fixed contrast we observe that Brinkman model is insensitive to
the distribution of the inclusions, the number of iterations is almost
the same for three permeability distributions: open foams, inclusions,
and connected inclusions.  However, in all permeability distributions
(except open foams) the iterations in both models grow with the
increase of the contrast. In the case (a) by increasing the contrast 100 times the
number of iterations increases by $50\%$ in Darcy model and at most
$100\%$ in Brinkman model. However, the numerical experiments 
shown on Table \ref{geom_set_reverse}  indicate that the iteration 
count is not very sensitive to the permeability distribution 
of case (b). These two experiments show that the preconditioner is quite sensitive
to the topology of the subregions with high and low permeability.

%
 \begin{figure}[tp]
        \centering
                  \begin{subfigure}[H]{0.2\textwidth}
       	                \includegraphics[height = 3.5cm, width=3.5cm]{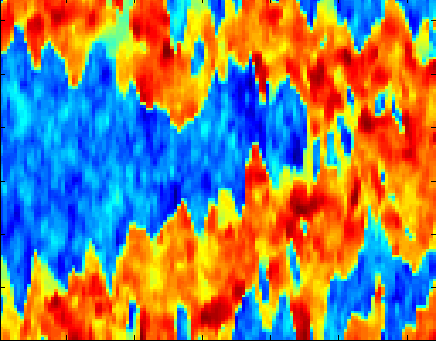}
                \caption{SPE10 Slice 44} 
        \end{subfigure}
        \hfil
          \begin{subfigure}[H]{0.2\textwidth}
       	                \includegraphics[height = 3.5cm, width=3.5cm]{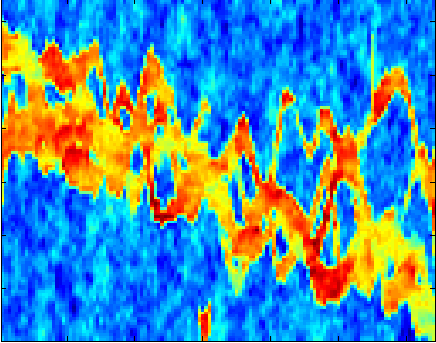}
                \caption{SPE10 Slice 49} 
        \end{subfigure}
        \hfil 
     \begin{subfigure}[H]{0.2\textwidth}
                    \includegraphics[height = 3.5cm, width=3.5cm]{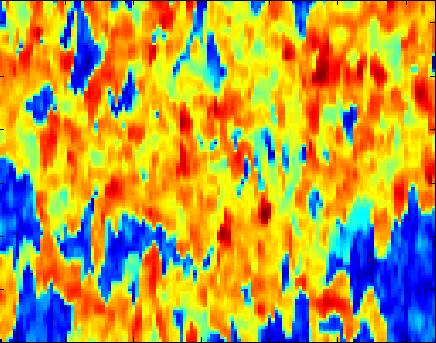}
                \caption{SPE10 Slice 54} 
        \end{subfigure}
        \hfil 
        \begin{subfigure}[H]{0.2\textwidth}
                    \includegraphics[height = 3.5cm, width=3.5cm]{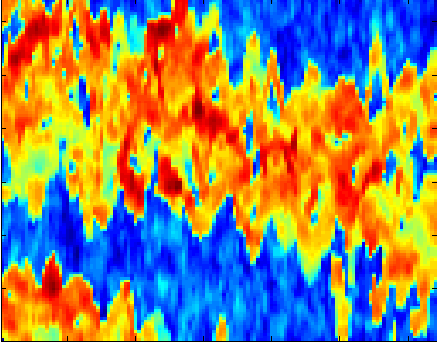}
                \caption{SPE10 Slice 74} 
        \end{subfigure}
              \caption{Logarithmic plots of the permeability in horizontal cross sections 44, 49, 54, and 74 of 
three dimensional  SPE10 benchmark \cite{SPE10_project}, Figure~\ref{fig:introduction}(B)   }\label{fig:set2}
\end{figure}  

\subsection{Tests with permeability of SPE10 benchmark}
The benchmark SPE10 \cite{ChristieBlunt01} 
provides permeability data on a 3-D sample represented on $60 \times 220 \times 85$ cells. 
There are 85 distinct layers within two general categories. 
The top 35 layers represent a prograding near shore environment 
with relatively low porosity and the bottom 50 layers represent a fluvial fan with channels.
The data set was rescaled on a $ 128 \times 128 \times 128$-mesh, presented in Figure~\ref{fig:introduction} (B).
Using $RT_0$ elements, namely the MAC scheme on this mesh, we obtain a
linear system of 15 million unknowns.  In order to run additional
examples, we have generated data for 2-D permeabilities, shown in
Figure~\ref{fig:set2}, that are horizontal slices of the SPE10
benchmark. These data sets were used to set up computational tests  for 2-D Darcy
and Brinkman models.  The permeability is plotted in logarithmic
scale, the contrast is almost $10^8$.  One can observe diverse
geological features in these slices: Slice 44, Figure~\ref{fig:set2}
(A), contains two distinct highly permeable channels starting from the
left that join together on the right; Slice 49 and 74,
Figure~\ref{fig:set2} (B) and (D), show a single highly permeable
channel (in red) with some small branches (or impermeable inclusions
within it), respectively; Slice 54, Figure~\ref{fig:set2} (C),
represent highly permeable region, which is almost dominated by the
red and yellow regions compared to Slices 44, 49, and 74.

\subsubsection{Performance of the iterative method for 2-D
  Darcy and Brinkman models}
Here we present iteration counts of the GMRES method applied to
Darcy and Brinkman models for permeability generated on horizontal
slices 44, 49, 54 and 74 of the SPE10 benchmark shown on
Figure~\ref{fig:set2}.

\begin{table}[tp]
\center
    \begin{tabular}{ | c | c c | c c | c c | c c | }
    \hline
     & \multicolumn{2}{c|}{Slice 44 }
     & \multicolumn{2}{c|}{Slice 49 }
     & \multicolumn{2}{c|}{Slice 54 }
     & \multicolumn{2}{c|}{Slice 74 }
    \\\hline
    Mesh size &$RT_0$&$RT_1$&$RT_0$&$RT_1$&$RT_0$&$RT_1$&$RT_0$&$RT_1$\\ \hline
  1/128&37 (36) & 42 (41) & 39 (36) &46 (42) &33 (32) &37 (35) &38 (36) &44 (43)
	  \\   \hline
 1/256&31 (31) &37 (37) &34 (32) &41 (38) &28 (28) &31 (30) &32 (31) &39 (38)
 \\ \hline
 1/512&27 (28) &33 (33) &31 (29) &37 (34) &24 (24) &27 (27) &27 (29) &35 (34)
 \\ \hline
    \end{tabular}
\vspace{2mm}
    \caption{Iteration count for Darcy and Brinkman (in parentheses) models with data of SPE10 benchmark 
(see Figure~\ref{fig:set2})  with RT$_0$   and RT$_1$ finite elements}
\label{darcy3}
 \end{table}

In Table~\ref{darcy3}, we present the number of iterations for  Darcy and Brinkman models
using $RT_0$ and $RT_1$ finite elements. 
For a given permeability distribution in a column, 
we give the number of GMRES iterations for mesh step-size varying from
$1/128$ to $1/512$ and for the permeability presented on Figure
~\ref{fig:set2}.  Clearly, the method is robust with respect to
mesh-size.  Moreover, analyzing the results in
Table 
\ref{darcy3} we see a reduction in the number of iterations for finer
meshes. A  similar observation can be made for the results of
Table~\ref{geom_set_1} as well.  A possible reason for such decrease
is that the preconditioner becomes better (as we refine the mesh) due
to the upscaling properties of the method, since the data is given on
fixed mesh $128 \times 128$).


We compared the iteration counts of our method with 
the recently proposed  in~\cite[formula 4.3]{Kraus_2015} block preconditioner of the Darcy system. In that paper
the block corresponding to the weighted $H^{\text{div}}$-norm 
(namely $(\kappa \qbu, \qbu) + (\nabla \cdot \qbu, \nabla \cdot \qbu)$) is handled by special 
multilevel method using additive Schur complement approximation.  
We compare the performance of these two methods on the permeability field generated by slice 44 of SPE10
(see, Figure~\ref{fig:set2}). 
The total number of  iterations 
in the minimal residual method in  \cite{Kraus_2015} is slightly less than the iterations shown on 
Table~\ref{darcy3}. However, within each iteration the method,
to invert the block corresponding to the weighted $H^{\text{div}}$-norm,
the method in  \cite[Tables 12 and 13]{Kraus_2015} uses in average
6 inner multilevel iterations.  
Our smoother is more expensive than the smoother in  \cite{Kraus_2015}, but the total number
of MG steps is much smaller. In summary, these two methods show similar performance and robustness
with respect to both step size and contrast. We note that the numerical tests show that our method performs 
efficiently on meshes with $2048 \times 2048$ (for RT$_0$) or $512 \times 512 $ (for RT$_1$)
for both Darcy and Brinkman models.
%

\subsubsection{Comparison of the solution of Darcy and Brinkman
  models}
\begin{figure}[tp]
          \centering
   \begin{subfigure}[H]{0.3\textwidth}
                \includegraphics[width=0.95\textwidth]
                {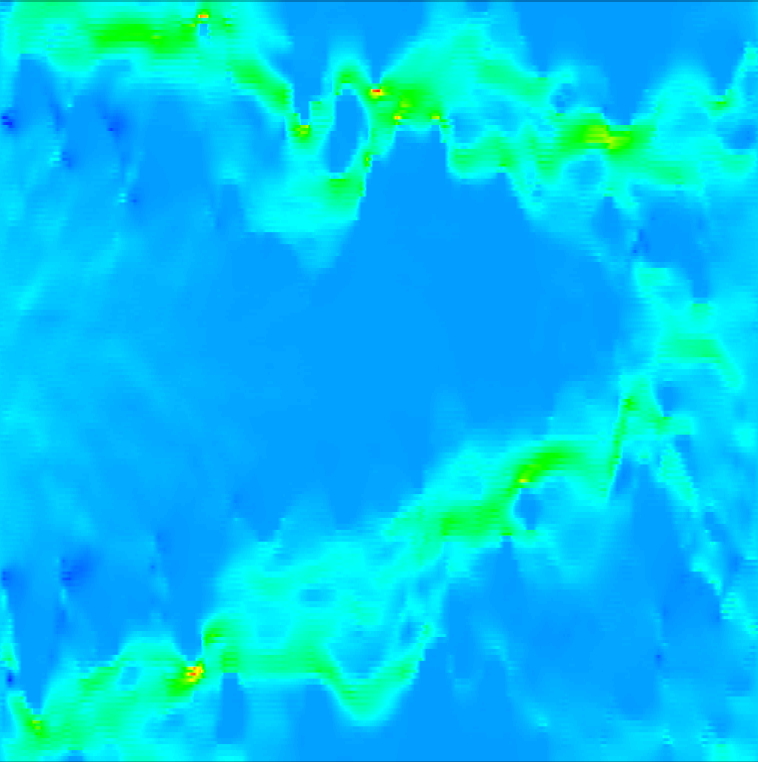}
                \caption{Darcy: $u_1$}  \label{u1-44} 
        \end{subfigure}
            \begin{subfigure}[H]{0.3\textwidth}
                \includegraphics[width=0.95\textwidth] 
                {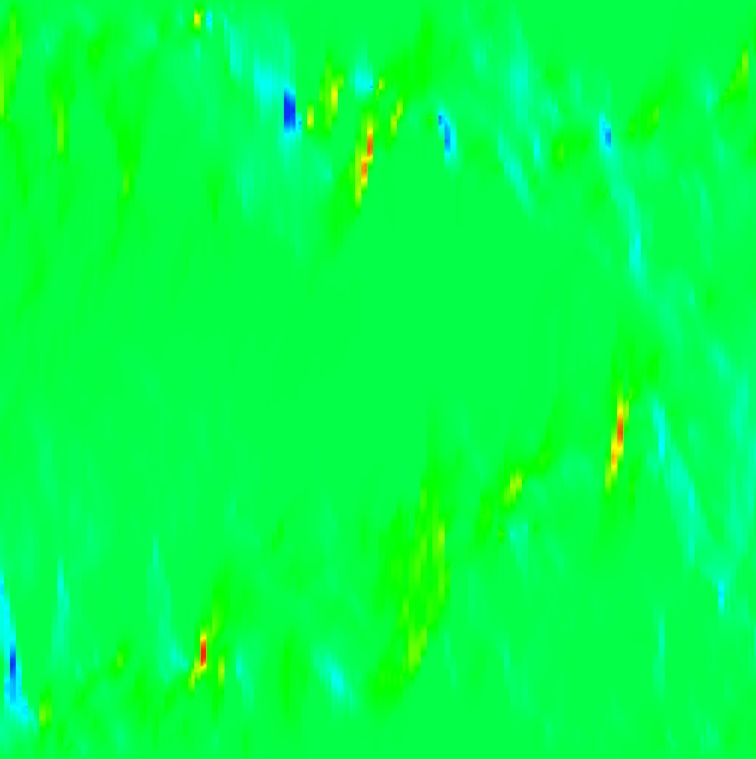}
                \caption{Darcy: $u_2$} \label{u2-44}
        \end{subfigure}
                 \begin{subfigure}[H]{0.3\textwidth}
                \includegraphics[width=0.95\textwidth] 
                {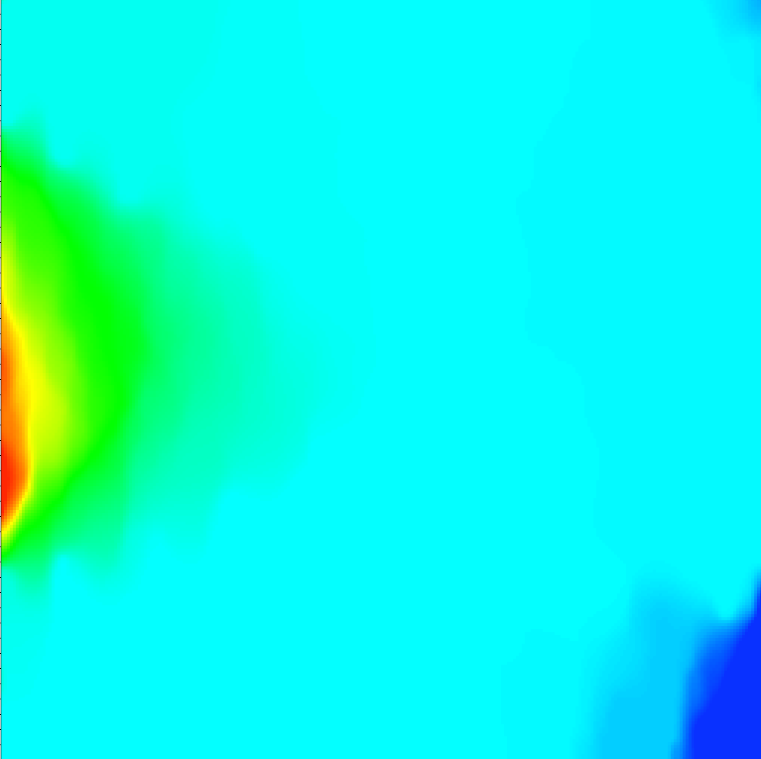}
                \caption{Darcy: $p$}   \label{p-44}
        \end{subfigure}        \\
        \vspace{2mm}
   \begin{subfigure}[H]{0.3\textwidth}
                \includegraphics[width=0.95\textwidth]{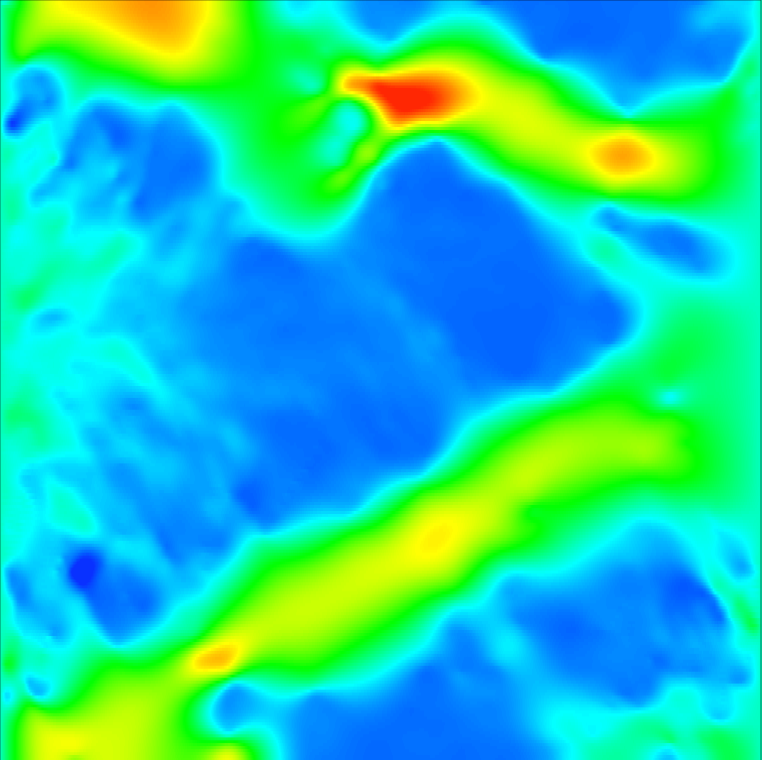}
                \caption{Brinkman: $u_1$}   \label{u1-44B}
        \end{subfigure}
        ~
            \begin{subfigure}[H]{0.3\textwidth}
                \includegraphics[width=0.95\textwidth]{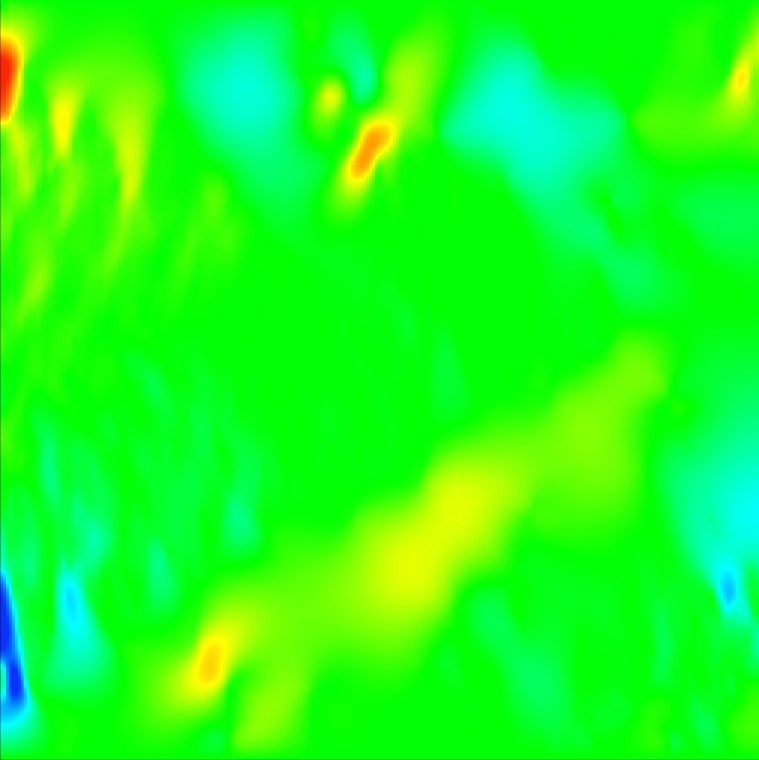}
                \caption{ Brinkman: $u_2$} \label{u2-44B}
        \end{subfigure}
                        ~
                 \begin{subfigure}[H]{0.3\textwidth}
                \includegraphics[width=0.95\textwidth]{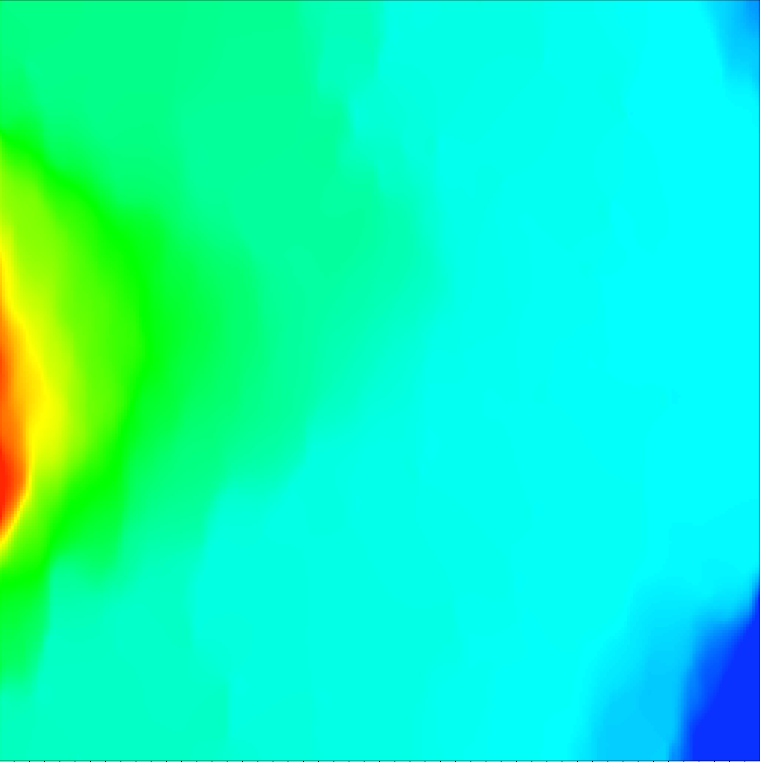}
                \caption{Brinkman:  $p$}   \label{p-44B}
        \end{subfigure}        
               \caption{
                Solution for permeability data of slice 44 in the same color scale} 
 \label{fig:spe10 d44}
\end{figure}

\begin{figure}[h!]
          \centering
%
        ~
   \begin{subfigure}[H]{0.3\textwidth}
                \includegraphics[width=0.95\textwidth] 
                {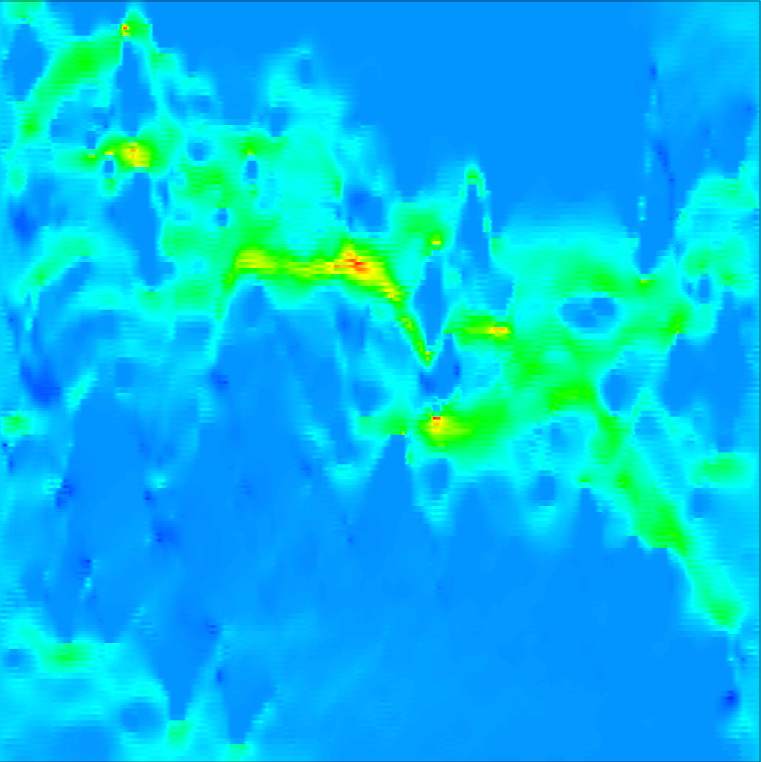}
                \caption{Darcy:  $u_1$}   \label{aa}
        \end{subfigure}
        ~
            \begin{subfigure}[H]{0.3\textwidth}
                \includegraphics[width=0.95\textwidth] 
                {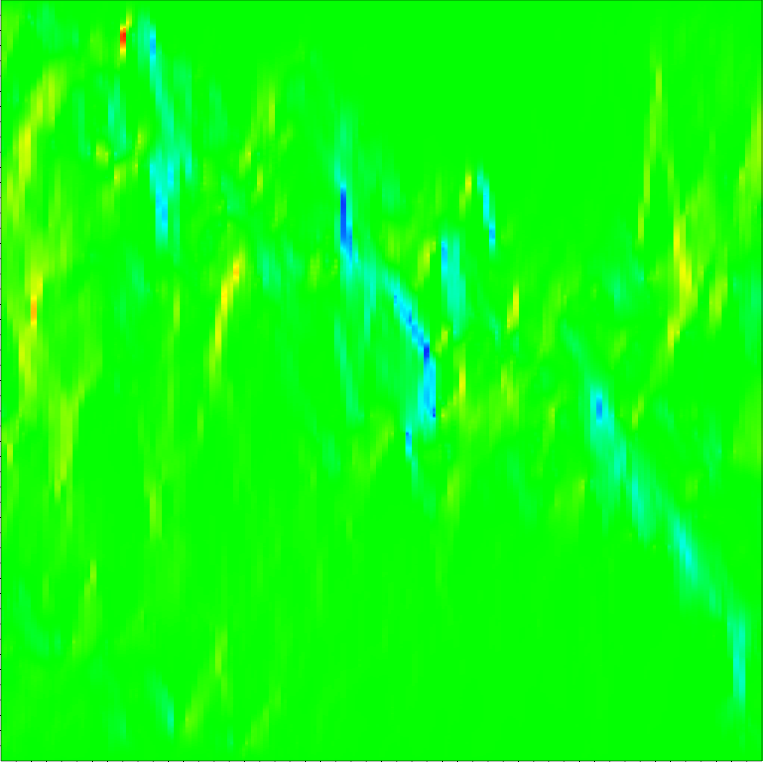}
                \caption{ Darcy:  $u_2$}
        \end{subfigure}
                        ~
                 \begin{subfigure}[H]{0.3\textwidth}
                \includegraphics[width=0.95\textwidth] 
                {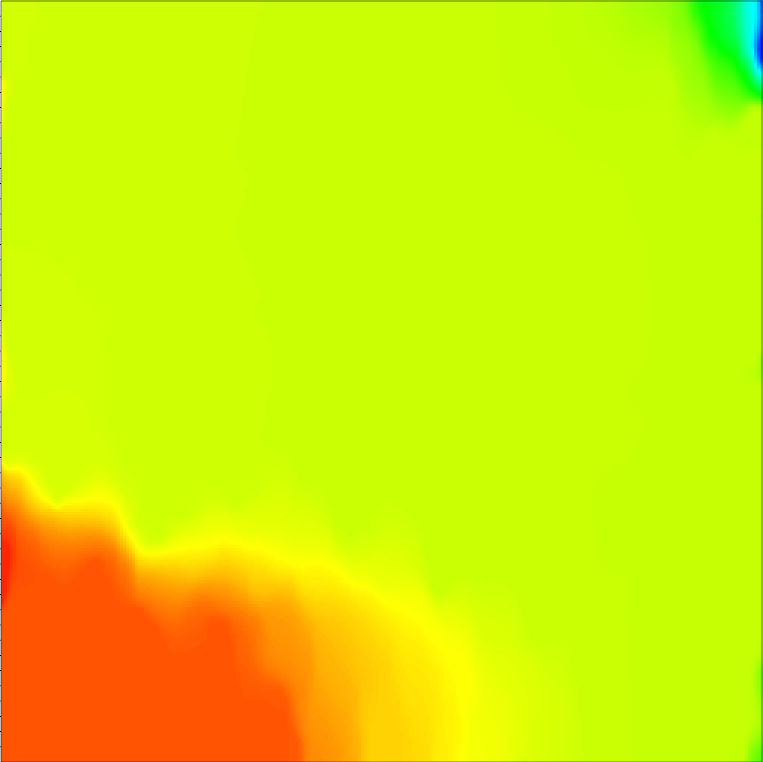}
                \caption{Darcy:  $ p$}   
        \end{subfigure}        
   \begin{subfigure}[H]{0.3\textwidth}
                \includegraphics[width=0.95\textwidth] 
                {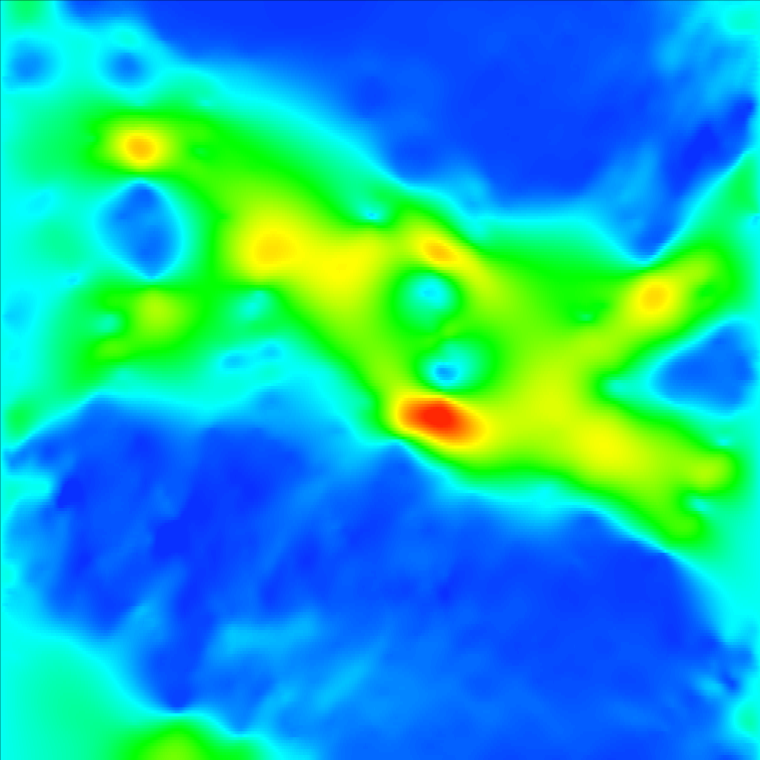}
                \caption{Brinkman:  $u_1$}   
        \end{subfigure}
        ~
            \begin{subfigure}[H]{0.3\textwidth}
                \includegraphics[width=0.95\textwidth] 
                {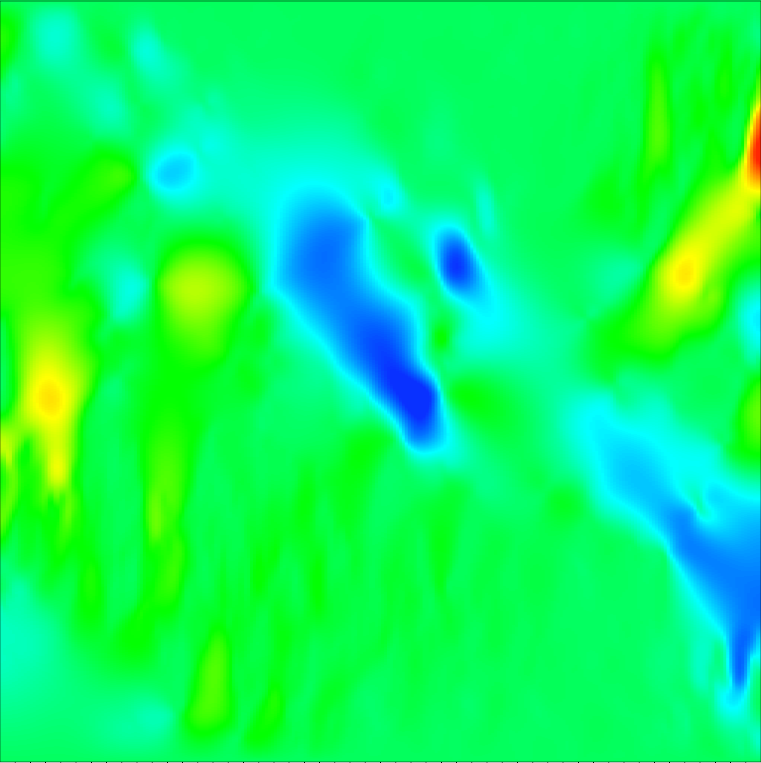}
                \caption{Brinkman: $u_2$}
        \end{subfigure}
                        ~
                 \begin{subfigure}[H]{0.3\textwidth}
                \includegraphics[width=0.95\textwidth] 
                {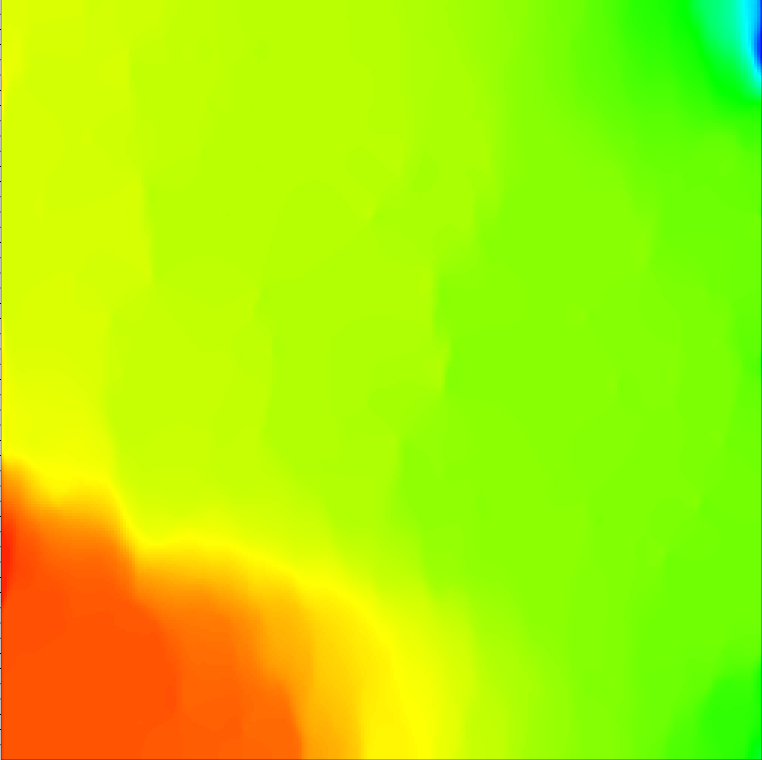}
                \caption{Brinkman: $ p$}   
        \end{subfigure}        
               \caption{
 Solution for permeability data of slice 74} 
 \label{fig:spe10 d74}
\end{figure}

Next, we present solutions produced by the two different models for
the permeability in Figures~\ref{fig:spe10 d44} and~\ref{fig:spe10 d74}, respectively.
Since the solutions are of low regularity we have
set all experiments for $RT_0$ finite elements.  The solutions of
Darcy and Brinkman models are plotted in the same scale so we can make
a visual comparison.  One can observe that the pressure obtained from
Darcy and Brinkman models in both slices, 44 and 49, do not differ
significantly. The most substantial differences are in the
velocity. 
One observes from this pictures that, as expected, the Brinkman model
is more diffusive.

\subsubsection{Computations with 3D SPE10 benchmark}
We now present and discuss briefly numerical results for the three
dimensional SPE10 benchmark, see, Figure~\ref{fig:introduction} (B).
We 
set $\mu = 10^{-2}$ and use Dirichlet boundary condition
${\qbg} =(1,0,0)$ and homogeneous right hand side, $ {\qbf}=\q0$.

In Table~\ref{brinkman3d}, we show the number of multigrid  iterations for two
different sets of boundary values and with three different
stopping criteria. One can observe that the method is quite robust with respect to different stopping criteria and boundary values. 

\begin{figure}[tp]
          \centering
    ~ 
                \begin{subfigure}[H]{0.22\textwidth}
                \includegraphics[height = 4cm, width=4.5cm]{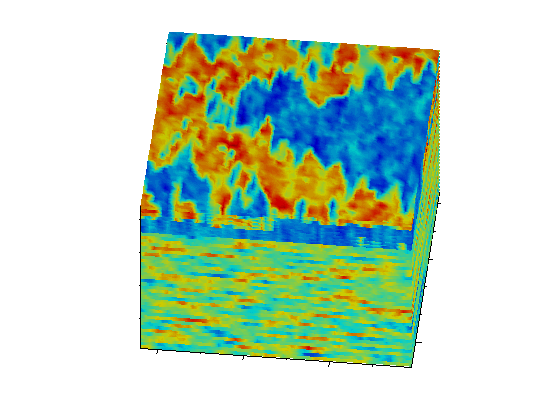}
                \caption{Slice 44: $\kappa(x)$} 
        \end{subfigure}
                ~
                  \begin{subfigure}[H]{0.22\textwidth}
                \includegraphics[height = 3.8cm, width=3.8cm]{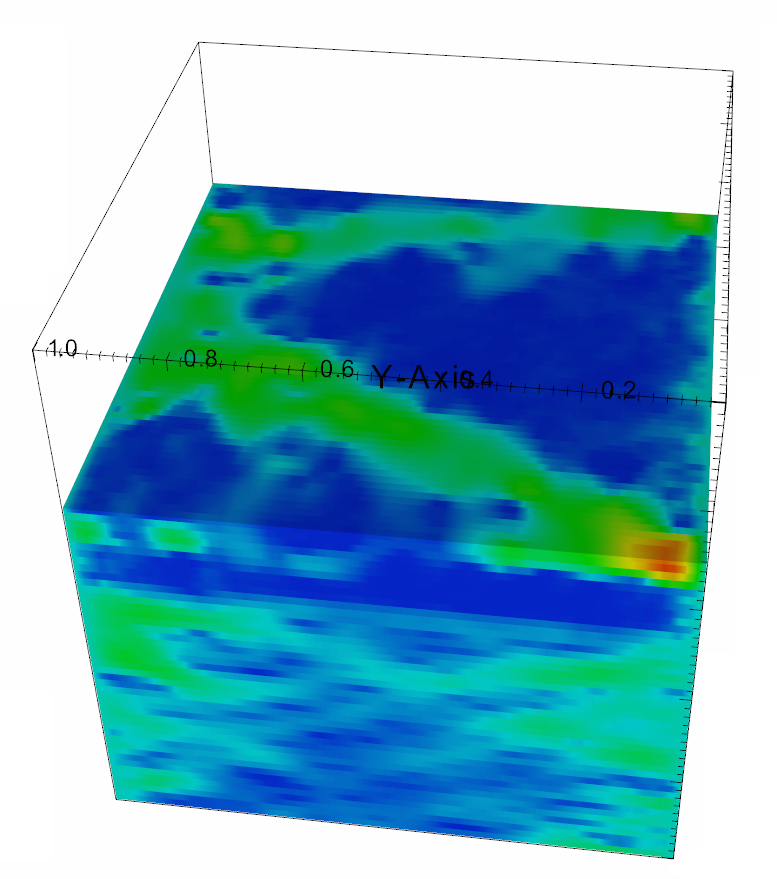}
                \caption{Slice 44: $u_1$} 
        \end{subfigure}
   ~
                \begin{subfigure}[H]{0.22\textwidth}
                \includegraphics[height = 4cm, width=4.5cm]{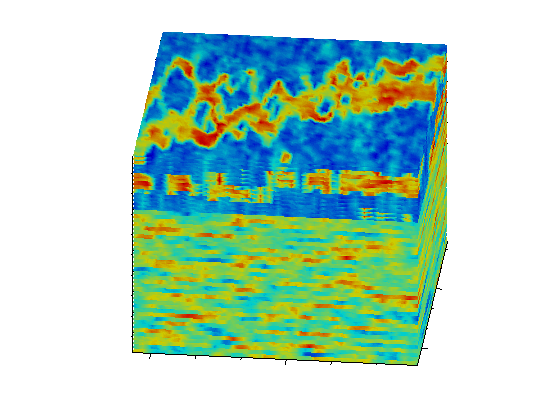}
                \caption{Slice 49:   $\kappa(x)$} 
        \end{subfigure}
                ~
                  \begin{subfigure}[H]{0.22\textwidth}
                \includegraphics[height = 3.8cm, width=3.8cm]{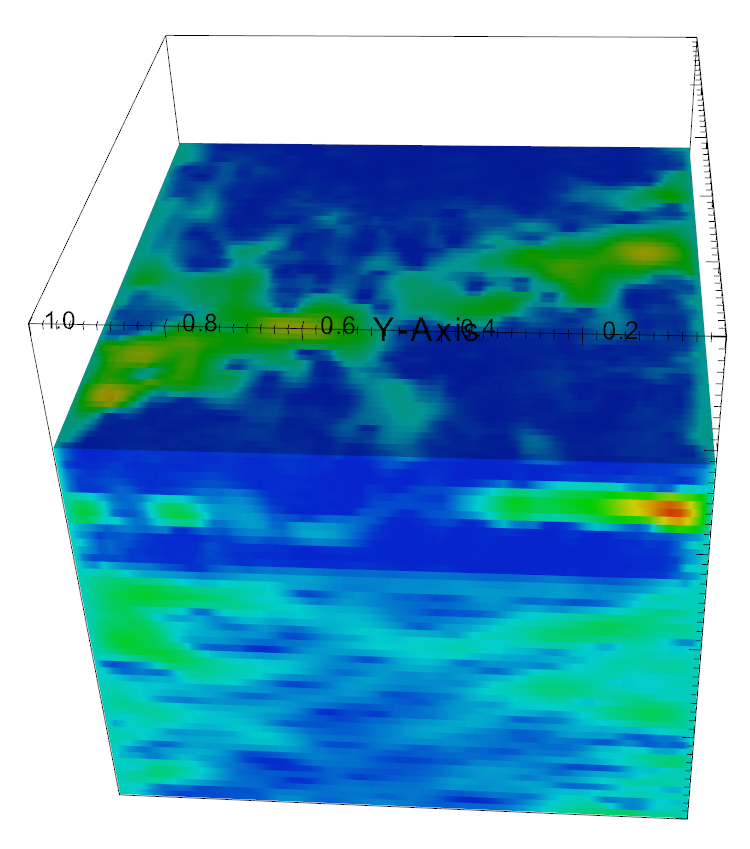}
                \caption{Slice 49: $u_1$} 
        \end{subfigure}
           \caption{Brinkman -  SPE10 with boundary data ${\qbg}=(1,0,0)$ } \label{yz-2}
\end{figure}

Finally, in Figure~\ref{yz-2} we show the $u_1$-component of the
solution in the horizontal planes corresponding to slices 44 and 49
obtained from 3-D Brinkman model solved with SPE10 benchmark
permeability. 

\begin{table}[tp]
\center
    \begin{tabular}{ | c | c c | c c | c c |}
    \hline
     & \multicolumn{2}{c|}{Residual $10^{-6}$}
        & \multicolumn{2}{c|}{Residual $10^{-8}$}
      & \multicolumn{2}{c|}{Residual $10^{-10}$}
    \\\hline
    Mesh size &$BC1$&$BC2$&$BC1$&$BC2$&$BC1$&$BC2$\\ \hline
  	 1/64&40&41&55&56&81&82
 \\ \hline
 1/128&31&31&46&46&73&73
 \\ \hline
    \end{tabular}
\vspace{2mm}
    \caption{Iteration count for the 3-D Brinkman model for SPE10 with boundary data -- 
 BC1:  ${\qbg}=(1,0,0)$ and BC2:  ${\qbg}=(0,1,0)$ }
\label{brinkman3d}
 \end{table}

\section{Conclusions}\label{sec:conclusions}


We studied the $H^{\text{div}}$-conforming discontinuous Galerkin
  discretization for  Darcy and Brinkman equations and presented a
  geometric multigrid method with a patch-based domain decomposition
  smoother for the corresponding algebraic system. Our numerical results show uniform
  contraction independent of the mesh size. In addition, we
  experimentally verify the robustness and efficiency of our method
  with respect to Darcy and Brinkman problems on two and three
  dimensional high contrast and high frequency distributions of the permeability. 
 Nevertheless, iteration counts go up into the 80s when we increase the contrast to $10^8$. We point
  out though,  that the geometric multigrid method with patch-based
  smoother is highly parallelizable.  We plan to study further some 
  algorithmic improvements like better choice of the coarse problem and
  parameter dependent projections.

\subsection*{Acknowledgments}

All computations were performed using the open source finite element
library deal.II~\cite{BangerthHartmannKanschat07,dealII83}.

\bibliographystyle{abbrv} 
\bibliography{references}

\end{document}